\documentclass[11pt,letterpaper]{article}
\usepackage{amsfonts, amsmath, amssymb, amscd, amsthm, graphicx}

\newtheorem{thm}{Theorem}[section]
\newtheorem{cor}[thm]{Corollary}
\newtheorem{lem}[thm]{Lemma}

\theoremstyle{definition}
\newtheorem{defn}[thm]{Definition}
\theoremstyle{remark}

\renewcommand{\phi }{{\rm\bf Lab\, }}
\newcommand{\e }{\varepsilon }
\newcommand{\G }{\Gamma (G, X\cup \mathcal H)}
\newcommand{\Ga }{\Gamma (G, \mathcal A)}
\newcommand{\dxh }{dist_{X\cup\mathcal H}}

\newcommand{\Hl }{\{ H_\lambda \} _{ \lambda \in \Lambda  }}

\newcommand{\C }{C(\e , \mu , \lambda , c, \rho )}

\begin{document}

\title{Small cancellations over relatively hyperbolic groups and embedding theorems}
\author{D. V. Osin \thanks{This work has been partially supported by the RFBR Grants $\sharp $ 02-01-00892, $\sharp $ 03-01-06555, and the NSF grant DMS-0605093}}
\date{}%

\maketitle

\begin{abstract}
We generalize the small cancellation theory
over ordinary hyperbolic groups to relatively hyperbolic settings.
This generalization is then used to prove various embedding theorem for
countable groups. For instance, we show that any countable torsion
free group can be embedded into a finitely generated group with
exactly 2 conjugacy classes. In particular, this gives the
affirmative answer to the well--known question of the existence of
a finitely generated group $G$ other than $\mathbb Z/2\mathbb Z$
such that all nontrivial elements of $G$ are conjugate.
\end{abstract}

%\tableofcontents

\parskip=.5mm
%%%%%%%%%%%%%%%%%%%%%%%%%%%%%%%%%%%%%%%%%%%%%%%%%%%%%%%%%%%%%

\section{Introduction}

%%%%%%%%%%%%%%%%%%%%%%%%%%%%%%%%%%%%%%%%%%%%%%%%%%%%%%%%%%%%%

Originally the notion of relative hyperbolicity was proposed by
Gromov \cite{Gro} in order to generalize various examples of
algebraic and geometric nature such as Kleinian groups,
fundamental groups of hyperbolic manifolds of pinched negative
curvature, small cancellation quotients of free products, etc. It
has been extensively studied in the last several years from
different points of view. The main aim of this paper is to
generalize the small cancellation theory over hyperbolic groups
developed by Olshanskii \cite{Ols2} to relatively hyperbolic
settings. Our approach is based on author's papers \cite{RHG,ESBG,RDF}, where the necessary background is provided. In the present paper we apply small cancellations over relatively hyperbolic groups to prove embedding theorems for countable groups. Further applications of our methods can be found in \cite{ABJ,SQ,BO,Min,Nor,Fac}.

In the paper \cite{HNN}, Higman, B.H. Neumann, and H. Neumann
proved that any countable group $G$ can be embedded into a
countable group $B$ such that every two elements of the same order
are conjugate in $B$. We notice that the group $B$ in \cite{HNN}
is constructed as a union of infinite number of subsequent
HNN--extensions and thus $B$ is never finitely generated. On the
other hand, any countable group can be embedded into a
$2$--generated group \cite{HNN}. Our first theorem is a natural
generalization of both of these results. For a group $G$, we
denote by $\pi (G)$ the set of all finite orders of elements of
$G$.

\begin{thm}\label{Conj}
Any countable group $G$ can be embedded into a $2$--generated
group $C$ such that any two elements of the same order are
conjugate in $C$ and $\pi (G)=\pi (C)$.
\end{thm}

\begin{cor}
Any countable torsion--free group can be embedded into a
(torsion--free) $2$--generated group with exactly $2$ conjugacy
classes.
\end{cor}

Since the number of finitely generated subgroups in any
2--generated groups is at most countable and the number of all
torsion--free finitely generated groups is uncountable, we have

\begin{cor}\label{2klassa}
There exists an uncountable set of pairwise non--isomorphic
torsion--free $2$--generated groups with exactly $2$ conjugacy
classes.
\end{cor}

We note that the question of the existence of {\it any}
finitely generated group with exactly $2$ conjugacy classes other
than $\mathbb Z/2\mathbb Z$ was open until now. It can be found,
for example, in \cite[Problem 9.10]{Kou} or in \cite[Problem
FP20]{Prob}. (Positive solution has been announced by Ivanov in 1989 \cite{Iva,IO0} but the complete proof has never been published.) Corollary \ref{2klassa} provides the first examples of such groups. Starting with the group $G=\mathbb Z/p^{n-2}\mathbb Z \times H$ for $n\ge 3$, where $p$ is a prime number and $H$ is a torsion--free group, we can generalize the previous result.

\begin{cor}\label{n}
For any $n\in \mathbb N$, $n\ge 2$, there is an uncountable set of
pairwise non--isomorphic finitely generated groups with exactly
$n$ conjugacy classes.
\end{cor}

For large enough prime numbers $n$, the first examples of finitely
generated infinite periodic groups with exactly $n$ conjugacy
classes were constructed by Ivanov (see \cite[Theorem 41.2]{Ols-book}) as limits
of hyperbolic groups (although hyperbolicity was not used
explicitly). Here we say that $G$ is a {\it limit of hyperbolic
groups} if there exists a finitely generated free group $F$ and a
series of normal subgroups $N_1\lhd N_2\lhd \ldots $ of $F$ such
that $G\cong F/N$ for $N=\bigcup\limits_{i=1}^\infty N_i$ and each
of the groups $F/N_i$, $i=1, 2, \ldots $ is hyperbolic. In
contrast it is impossible to construct a finitely generated group
other than $\mathbb Z/2\mathbb Z$ with exactly $2$ conjugacy
classes in this way.

Indeed suppose that a finitely generated group $G$ has exactly $2$
conjugacy classes. If $G$ is not torsion--free, then $G$ is a
group of exponent $p$ for some prime $p$ as the orders of all
nontrivial elements of $G$ are equal. If $p=2$, $G$ is abelian and
hence is isomorphic to $\mathbb Z/2\mathbb Z$. In case $p>2$,
there exist non--trivial elements $g,t\in G $ such that
\begin{equation} \label{BS}
t^{-1}gt=g^2.
\end{equation}
The equality $g^{2^p-1}=t^{-p}gt^pg^{-1}=gg^{-1}=1$ implies
$2^p-1\equiv 0(mod\; p)$. However, by the Fermat Little Theorem,
we have $2^p-2\equiv 0(mod\; p)$, which contradicts the previous
equality. Assume now that $G$ is torsion--free. If $G$ is a limit
of hyperbolic groups $F/N_i$, $i=1,2,\ldots $, then for some $i$
large enough, there are elements $g,t\in F/N_i$ of infinite order
that satisfy (\ref{BS}). This leads to a contradiction again since
the equality of type (\ref{BS}) is impossible in a hyperbolic
group if the order of $g$ is infinite \cite{GH,Gro}.

Another theorem from \cite{HNN} states that any countable group
$G$ can be embedded into a countable divisible group $D$. We
recall that a group $D$ is said to be {\it divisible} if for every
element $d\in D$ and every positive integer $n$, the equation
$x^n=d$ has a solution in $D$. A natural example of a divisible
groups is $\mathbb Q$. The question of
the existence of a finitely generated divisible group was open
during a long time. The first examples of such a type were
constructed by Guba \cite{Guba} (see also \cite{Ols-book}).

Later Mikhajlovskii and Olshanskii \cite{MO} constructed a more
general example of a finitely generated {\it verbally complete
group}, that is a group $W$ such that for every nontrivial freely reduced word $w(x_i)$ in
the alphabet $x_1^{\pm 1} ,x_2^{\pm 1},\ldots $ and every $v\in
W$, the equation $w(x_i)=v$ has a solution in $W$. That is, there
are elements $v_1, v_2, \ldots \in W$ such that $w(v_i)=v$ in $W$,
where $w(v_i)$ is the word obtained from $w(x_i)$ by substituting
$v_i$ for $x_i$, $i=1,2 \ldots $. Comparing these results one may
ask whether any countable group can be embedded into a finitely
generated divisible (or verbally complete) group. The next theorem
provides the affirmative answer.

\begin{thm} \label{WAC}
Any countable group $H$ can be embedded into a $2$--generated
verbally complete group $W$. Moreover, if $H$ is torsion--free,
then $W$ can be chosen to be torsion--free.
\end{thm}

Note that the condition $\pi (G)=\pi (W)$ can not be ensured in
Theorem \ref{WAC}. Indeed, it is easy to show that if a divisible
group $W$ contains a nontrivial element of finite order, then $\pi
(W)=\mathbb N$. As above, we obtain

\begin{cor} \label{WAC1}
There exists an uncountable set of pairwise non--isomorphic
$2$--generated verbally complete groups.
\end{cor}

{\bf Acknowledgment.} The author is grateful to A. Minasyan and A. Olshanskii for helpful discussions and to the anonymous referee for careful reading of the manuscript and useful remarks.

%%%%%%%%%%%%%%%%%%%%%%%%%%%%%%%%%%%%%%%%%%%%%%%%%%%%%%%%%%%%%%%%%

\section{Outline of the method}

%%%%%%%%%%%%%%%%%%%%%%%%%%%%%%%%%%%%%%%%%%%%%%%%%%%%%%%%%%%%%%%%%

In this section we give the proofs of Theorem \ref{WAC} and
Theorem \ref{Conj} modulo technical results which are obtained in
Sections 4-8. We assume the reader to be familiar with the notion
of a relatively hyperbolic group and refer to the next section for
precise definitions.

Let $G$ be a group that is hyperbolic relative to a collection of
subgroups $\Hl $. We divide the set of all elements of $G$ into
two subsets as follows. An element $g\in G$ is said to be {\it
parabolic} if $g$ is conjugate to an element of $H_\lambda $ for
some $\lambda \in \Lambda $. Otherwise $g$ is said to be {\it
hyperbolic}. Recall also that a group is {\it elementary } if it
contains a cyclic subgroup of finite index. The following result
concerning maximal elementary subgroups is proved in \cite[Theorem 4.3, Corollary 1.7]{ESBG}.

\begin{thm}\label{E(g)}
Let $G$ be a group hyperbolic relative to a collection of
subgroups $\Hl $, $g$ a hyperbolic element of infinite order of
$G$. Then the following conditions hold.
\begin{enumerate}
\item The element $g$ is contained in a unique maximal elementary
subgroup $E_G(g)$ of $G$, where $$E_G(g)=\{ f\in G\; :\;
f^{-1}g^nf=g^{\pm n}\; {\rm for \; some\; } n\in \mathbb N\} .$$

\item The group $G$ is hyperbolic relative to the collection
$\Hl\cup \{ E_G(g)\} $.
\end{enumerate}
\end{thm}

Given a subgroup $H\le G$, we denote by $H^0$ the set of all
hyperbolic elements of infinite order in $H$. Recall also that two
elements $f,g\in G^0$ are said to be {\it commensurable} (in G) if
$f^k$ is conjugate to $g^l$ in $G$ for some non--zero $k,l$.

\begin{defn} \label{suit} A subgroup $H\le G$ is called
{\it suitable}, if there exist two non--commensurable elements
$f_1, f_2\in H^0$ such that $E_G(f_1)\cap E_G(f_2)=1$.
\end{defn}

The next lemma is proved in Section 8.

\begin{lem}\label{non-com}
Let $G$ be a group hyperbolic relative to a collection of
subgroups $\Hl $, $H$ a suitable subgroup of $G$. Then there exist
infinitely many pairwise non--commensurable (in G) elements $h_1,
h_2, \ldots \in H^0$ such that for all $i=1,2, \ldots $,
$E_G(h_i)=\langle h_i\rangle $. In particular, $E_G(h_i)\cap
E_G(h_j)=\{ 1\} $ whenever $i\ne j$.
\end{lem}

Our main tool is the following theorem proved in Section 8. The
proof is based on a certain small cancellation techniques developed
in Sections 4-7.

\begin{thm}\label{glue}
Let $G$ be a group hyperbolic relative to a collection of
subgroups $\Hl $, $H$ a suitable subgroup of $G$, and $t_1, \ldots
, t_m$ arbitrary elements of $G$. Then there exists an epimorphism
$\eta \colon G\to \overline{G}$ such that:
\begin{enumerate}
\item The group $\overline{G}$ is hyperbolic relative to $\{ \eta
(H_\lambda ) \} _{\lambda \in \Lambda } $. \item For any $i=1,
\ldots , m$, we have $\eta (t_i)\in \eta (H)$. \item The
restriction of $\eta $ to $\bigcup\limits_{\lambda\in \Lambda
}H_\lambda $ is injective. \item $\eta (H)$ is a suitable subgroup
of $\overline{G}$. \item Every element of finite order in $\overline{G} $ is an image of an element of finite order in
$G$. In particular, if all hyperbolic elements of $G$ have
infinite order, then all hyperbolic elements of $\overline G$ have
infinite order.
\end{enumerate}
\end{thm}

The next theorem is proved in \cite{HNN}[Corollary 1.4]. For
finitely generated groups this result was first proved by Dahmani
in \cite{Dah}. It is worth noting that we use the theorem for
infinitely generated groups in this paper.

\begin{thm}\label{HNN0}
Suppose that a group $G$ is hyperbolic relative to a collection of
subgroups $\Hl \cup \{ K\} $ and for some $\nu \in \Lambda $,
there exists a monomorphism $\iota \colon K\to H_{\nu }$. Then the
HNN--extension
\begin{equation}\label{HNN-pres}
\widetilde G=\langle G,t\; |\; t^{-1}kt=\iota (k),\; k\in K\rangle
\end{equation}
is hyperbolic relative to $\Hl $.
\end{thm}

Theorems \ref{WAC} and \ref{Conj} can be obtained in a uniform way
from the following result.

\begin{thm}\label{0}
Suppose that $R$ is a countable group such that for any elementary
group $E$ satisfying the condition $\pi (E)\subseteq \pi (R)$,
there exists a subgroup of $R$ isomorphic to $E$. Then there is an
embedding of $R$ into a 2--generated group $S=S(R)$ such that any
element of $S$ is conjugate to an element of $R$ in $S$. In
particular, $\pi (S)=\pi (R)$.
\end{thm}

\begin{proof}
The desired group $S$ is constructed as an inductive limit of
relatively hyperbolic groups as follows. Let us set $$G(0)=R\ast
F(x,y),$$ where $F(x,y)$ is the free group of rank $2$ generated
by $x$ and $y$. We enumerate all elements of $$R=\{ 1=r_0, r_1,
r_2, \ldots \} $$ and $$G(0)=\{ 1=g_0, g_1, g_2, \ldots \} .$$

Suppose that for some $i\ge 0$, the group $G(i)$ has already been
constructed together with an epimorphism $\xi _i\colon G(0)\to
G(i)$. We use the same notation for elements $x,y, r_0, r_1,
\ldots , g_0, g_1, \ldots $ and their images under $\xi _i $ in
$G(i)$. Assume that $G(i)$ satisfies the following conditions. (It
is straightforward to check these conditions for $G(0)$ and the
identity map $\xi _0\colon G(0)\to G(0)$.)
\begin{enumerate}
\item[(i)] The restriction of $\xi _i$ to the subgroup $R$ is
injective. In what follows we identify $R$ with its image in
$G(i)$.

\item[(ii)] $G(i)$ is hyperbolic relative to $R$.

\item[(iii)] The elements $x$ and $y$ generate a suitable subgroup
of $G(i)$.

\item[(iv)] All hyperbolic elements of $G(i)$ have infinite order.
In particular, $\pi (G(i))=\pi (R)$.

\item[(v)] The elements $g_0, \ldots , g_i$ are parabolic in
$G(i)$.

\item[(vi)] In the group $G(i)$, the elements $r_0, \ldots , r_i$
are contained in the subgroup generated by $x$ and $y$.
\end{enumerate}
The group $G(i+1)$ is obtained from $G(i)$ in two steps.

{\bf Step 1.} Let us take the element $g_{i+1}$ and construct a
group $G(i+1/2)$ as follows. If $g_{i+1}$ is a parabolic element
of $G(i)$, we set $G(i+1/2)=G(i)$. If $g_{i+1}$ is hyperbolic, the
order of $g_{i+1}$ is infinite by (iv). Furthermore, since $\pi
(E_{G(i)}(g_{i+1}))\subseteq \pi (G(i))=\pi (R)$, there is a
monomorphism $\iota\colon E_{G(i)}(g_{i+1})\to R$. Then we take
the HNN--extension
$$G(i+1/2 )=\langle G(0), t\; | \; t^{-1}et=\iota (e), \, e\in
E_{G(i)}(g_{i+1})\rangle .$$

In both cases $G(i+1/2)$ is hyperbolic relative to $R$. Indeed
this is obvious in the first case and follows from the second
assertion of Theorem \ref{E(g)} and Theorem \ref{HNN0} in the
second one. Note also that all hyperbolic elements of $G(i+1/2)$
have infinite order. (In the second case this immediately follows
from the description of periodic elements in HNN--extensions
\cite[Ch. IV, Theorem 2.4]{LS}.)

{\bf Step 2.} First we wish to show that the subgroup generated by
$x$ and $y$ is suitable in $G(i+1/2)$. This is obvious in case
$g_{i+1}$ is parabolic in $G(i)$, so we consider the second case
only. Since $\langle x, y\rangle $ is suitable in $G(i)$ by (iii),
Lemma \ref{non-com} yields the existence of infinitely many
pairwise non--commensurable (in $G(i)$) hyperbolic elements
$h_j\in\langle x, y\rangle $ of infinite order, $j=1,2 \ldots $,
such that $E_{G(i)}(h_j)=\langle h_j\rangle $. At most one of
these elements is commensurable with $g_{i+1}$ in $G(i)$.
Therefore, there exist two non--commensurable in $G(i)$ hyperbolic
elements of infinite order, say $h_1,h_2\in \langle x, y\rangle $,
such that $h_j$ is not commensurable with $g_{i+1}$ in $G(i)$ for
$j=1,2$. In particular, $h_j$, $j=1,2$, is not conjugate to an
element of $E_{G(i)}(g_{i+1})$ as $\langle g_{i+1}\rangle $ has
finite index in $E_{G(i)}(g_{i+1})$. According to Britton's Lemma
on HNN--extensions \cite[Ch. 5, Sec.2]{LS}, this implies that
$h_1$ and $h_2$ are hyperbolic and non--commensurable in
$G(i+1/2)$. Furthermore, if for some $j=1,2$, $n\in \mathbb N$,
and $u\in G(i+1/2)$, we have $u^{-1}h_j^nu=h_j^{\pm n}$, then
$u\in G(i)$ by Britton's Lemma. Thus the explicit description of
maximal elementary subgroups from the first assertion of Theorem
\ref{E(g)} yields the equality $E_{G(i+1/2)}(h_j)=E_{G(i)}(h_j)$
for $j=1,2$. Finally since $E_{G(i)}(h_j)=\langle h_j\rangle $ and
$h_1$, $h_2$ are non--commensurable, we have
$$E_{G(i+1/2)}(h_1)\cap E_{G(i+1/2)}(h_2)=E_{G(i)}(h_1)\cap
E_{G(i)}(h_2)=\langle h_1\rangle \cap \langle h_2\rangle=\{ 1\}
.$$ By Definition \ref{suit} this means that the subgroup
generated by $x$ and $y$ is suitable in $G(i+1/2)$.

We now apply Theorem \ref{glue} to the group $G=G(i+1/2)$, the
subgroup $H=\langle x,y\rangle \le G(i+1/2)$, and the set of
elements $\{ t, r_{i+1}\} $. Let $G(i+1)= \overline{G}$, where
$\overline {G}$ is the quotient group provided by Theorem
\ref{glue}. Since $t$ becomes an element of $\langle x,y\rangle $
in $G(i+1)$, there is a naturally defined epimomorphism $\xi
_{i+1}\colon G(0)\to G(i+1)$. Using Theorem \ref{glue} it is
straightforward to check properties (i)--(vi) for $G(i+1)$. This
completes the inductive step.

Let $N_i$ denote the kernel of $\xi _i $. Observe that $N_1, N_2,
\ldots $ form an increasing normal series and set $S=G(0)/N$,
where $N=\bigcup_{i=1}^\infty N_i$. By (i) the subgroup $R$ is
embedded into $S$. Further it is easy to see that $S$ is
$2$--generated. Indeed, $G(0)$ is generated by $x,y, r_1, r_2,
\ldots $. Condition (vi) yields $ r_i\in \langle x,y \rangle $ in
$S$ for any $i\in \mathbb N$. Thus $S$ is generated by $x$ and
$y$. Finally let $s$ be an element of $S$. We take an arbitrary
preimage $g\in G(0)$ of $s$. Then the image of the element $g$
becomes parabolic at a certain step according to (v). Thus $s$ is
conjugate to an element of $R$ in $S$. The theorem is proved.
\end{proof}

It remains to derive Theorems \ref{WAC} and \ref{Conj}.

\begin{proof}[Proof of Theorem \ref{Conj}.]
Let $\mathcal E$ denote the free product of all elementary groups
$E$ (taken up to isomorphism) such that $\pi (E)\subseteq \pi
(G)$. We set $G^\ast =G\ast \mathcal E$. By a theorem from
\cite{HNN}, we can embed $G^\ast $ into an (infinitely generated)
group $R$ such that all elements of the same order are conjugate
in $R$ and
\begin{equation}\label{pi}
\pi (R)=\pi (G^\ast )=\pi (G).
\end{equation}

We now apply Theorem \ref{0} and embed the group $R$ into a
2-generated group $C=S(R)$ such that any element of $C$ is
conjugate to an element of $R$. As all elements of the same order
are  conjugate in $R$, this is so in $C$. The equality $\pi
(C)=\pi (G)$ follows from (\ref{pi}) as $\pi (C)=\pi (R)$ by
Theorem \ref{0}.
\end{proof}

\begin{proof}[Proof of Theorem \ref{WAC}.]
First note that any countable group $G$ can be embedded into an
infinitely generated countable verbally complete group $R$ in the
following way. (The idea comes from the proof of the
Higman--Neumann--Neumann theorem on embeddings into divisible
groups.) We denote by $F=F(a_1, a_2, \ldots )$ the free group with
basis $a_1, a_2, \ldots $. Let us enumerate the set of all pairs
$$\{ p_1, p_2, \ldots \} =\{ (v, g)\; :\; v\in
F\setminus \{ 1\} ,\, g\in G\setminus\{ 1\} \} .$$ Starting with
the group $G$ we first set $G^\ast =G$ if $G$ is torsion--free,
and $G^\ast =G\ast E_1\ast E_2 \ast \ldots $, where $\{ E_1, E_2,
\ldots \} $ is the set of all elementary groups (up to
isomorphism), otherwise. Further we construct a sequence of groups
$G^\ast=U_0\le U_1\le \ldots $ as follows. Suppose that for some $i\ge
0$, the group $U_i$ has already been constructed and take
$p_{i+1}=(v,g)$. There are two possibilities to consider.

1) The element $g$ has infinite order. Then we define $U_{i+1}$ to
be the free product of $U_i$ and $F$ with the amalgamated
subgroups $\langle g\rangle $ and $\langle v \rangle $.

2) The order of $g$ is $n< \infty $. It is well--known
\cite[Theorem 5.2, Ch. 4]{LS} that the order of the element $v$ in
the group $H=\langle a_1, a_2, \ldots \; |\; v^n=1\rangle $ equals
$n$. Thus the free product of $U_i$ and $H$ with amalgamated
subgroups $\langle g\rangle $ and $\langle v \rangle $ is
well--defined. We set $U_{i+1}=U_i\ast _{\langle g\rangle =\langle
v \rangle } H$.

Now let $U(G^\ast)=\bigcup\limits_{i=0}^\infty U_i $. Obviously $G^\ast$
embeds in $U(G^\ast)$, $U(G^\ast)$ is countable and torsion--free whenever
$G^\ast$ is torsion--free, and any equation of type $w(x_i)=g$, where
$w(x_i)$ is a word in the alphabet $x_1^{\pm 1}, x_2^{\pm 1},
\ldots $ and $g\in G$, has a solution in $U(G^\ast)$. Finally we
consider the sequence of groups $R_1\le R_2 \le \ldots $, where
$R_1=U(G^\ast)$ and $R_{i+1}=U(R_i)$, $i=1,2 \ldots $. Clearly the
group $R=\bigcup\limits_{i=0}^\infty R_i $ is countable, verbally
complete, torsion--free whenever $G$ is torsion--free, and
contains a copy of every elementary group $E$ such that $\pi
(E)\subseteq \pi (G)$. Let $W=S(R)$ be the group provided by
Theorem \ref{0}.

Consider an equations $w(x_i)=v$ for some $v\in W$. By Theorem
\ref{0}, there is an element $t\in W$ such that $t^{-1}vt\in R$.
Since $R$ is verbally complete, there is a solution $x_1=r_1$,
$x_2=r_2$, $\ldots $ to the equation $w(x_i)=t^{-1}vt $ in $R$.
Clearly $x_1=tr_1t^{-1}$, $x_2=tr_2t^{-1}$, $\ldots $ is a
solution to the equation $w(x_i)=v$.
\end{proof}

%%%%%%%%%%%%%%%%%%%%%%%%%%%%%%%%%%%%%%%%%%%%%%%%%%%%%%%%%%%%%%%%

\section{Preliminaries}

%%%%%%%%%%%%%%%%%%%%%%%%%%%%%%%%%%%%%%%%%%%%%%%%%%%%%%%%%%%%%%%%

{\bf Some conventions and notation.} We write $W\equiv V$ to express the letter--for--letter equality of
words $W$ and $V$ in some alphabet. If a word $W$ decomposes as $W\equiv V_1UV_2$, we call $V_1$ (respectively, $V_2$) a {\it prefix} (respectively, {\it suffix}) of $W$.
For elements $g$, $t$ of a group $G$, $g^t$
denotes the element $t^{-1}gt$. Recall that a subset $X$ of a
group $G$ is said to be {\it symmetric} if for any $x\in X$, we
have $x^{-1}\in X$. In this paper all generating sets of groups
under consideration are supposed to be symmetric.

All paths considered in this paper are combinatorial paths. Recall that a {\it combinatorial path} $p$ in a $CW$-complex is a sequence of edges (i.e., 1-dimensional cells) $e_1e_2\ldots e_k$, where $(e_i)_+=(e_{i+1})_-$. If edges of the complex are labeled, we define {\it the label of $p$} by $\phi (p)\equiv \phi (e_1)\phi(e_2)\ldots \phi (e_k),$ where $\phi (e_i)$ is the label of $e_i$.  We also denote by $p_-=(e_1)_-$
and $p_+=(e_k)_+$ the origin and the terminus of $p$ respectively.
The length $l(p)$ of $p$ is the number of edges of $p$.

\vspace{3mm}

\noindent {\bf Word metrics and Cayley graphs.} Let $G$ be a group
generated by a (symmetric) set $\mathcal A$. Recall that the {\it
Cayley graph} $\Ga $ of a group $G$ with respect to the set of
generators $\mathcal A$ is an oriented labeled 1--complex with
the vertex set $V(\Ga )=G$ and the edge set $E(\Ga )=G\times
\mathcal A$. An edge $e=(g,a)$ goes from the vertex $g$ to the
vertex $ga$ and has label $\phi (e)\equiv a$. As usual, we denote
the origin and the terminus of the edge $e$, i.e., the vertices
$g$ and $ga$, by $e_-$ and $e_+$ respectively.

Associated to $\mathcal A$ is the so--called {\it word metric} on
$G$. More precisely, the length $|g|_\mathcal A$ of an element
$g\in G$ is defined to be the length of a shortest word in
$\mathcal A$ representing $g$ in $G$. By abuse of notation, we also write
$|W|_\mathcal A$ to denote the lengths of the element of $G$ represented by a word $W$ in the alphabet $\mathcal A$. This is
to be distinguished from the lengths of the word $W$ itself, which is denoted by $\| W\| $.

The word metric on $G$
is defined by $dist_\mathcal A(f,g)=|f^{-1}g|_\mathcal A$. We also denote by
$dist _\mathcal A$ the natural extension of the word metric to the
Cayley graph $\Ga $.

\vspace{3mm}

\noindent {\bf Van Kampen Diagrams.} Recall that a {\it van Kampen
diagram} $\Delta $ over a presentation
\begin{equation}
G=\langle \mathcal A\; | \; \mathcal O\rangle \label{ZP}
\end{equation}
is a finite oriented connected 2--complex endowed with a labelling
function $\phi : E(\Delta )\to \mathcal A$, where $E(\Delta ) $
denotes the set of oriented edges of $\Delta $, such that $\phi
(e^{-1})\equiv (\phi (e))^{-1}$. Labels and lengths of paths are
defined as in the case of Cayley graphs. Given a cell $\Pi $ of
$\Delta $, we denote by $\partial \Pi$ the boundary of $\Pi $;
similarly, $\partial \Delta $ denotes the boundary of $\Delta $.
The labels of $\partial \Pi $ and $\partial \Delta $ are defined
up to a cyclic permutation. An additional requirement is that for
any cell $\Pi $ of $\Delta $, the boundary label $\phi (\partial
\Pi)$ is equal to a cyclic permutation of a word $P^{\pm 1}$,
where $P\in \mathcal O$.

\begin{figure}
\vspace{2mm}\hspace{15mm} \includegraphics{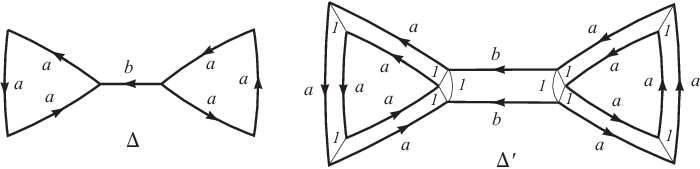}\\
\vspace{-2mm}
 \caption{$\Delta ^\prime$ is a $0$--refinement of the diagram $\Delta $
over the presentation $\langle a,b\; |\; a^3=1\rangle $.}
\label{figa}
\end{figure}

Sometimes it is convenient to use the
notion of the so-called {\it $0$--refinement}, which enables us to assume
that all diagrams a homeomorphic to a disc. Roughly speaking, making a $0$-refinement of a diagram $\Delta $
just means replacing every edge $e\in E(\Delta )$ with a rectangle labeled $\phi (e) 1\phi(e)^{-1} 1$, and replacing every vertex of $\Delta $ with a polygon (2-cell) labeled by $11\cdots 1$. Here $1$ means the empty word. The rectangles are then attached to the polygons along edges labeled by $1$ (see Fig. \ref{figa}). The lengths of edges labelled by $1$ is supposed to be $0$. This notion is quite standard and we refer the reader to \cite[Ch. 4]{Ols-book} for details.

The van Kampen Lemma states that a word $W$ over the alphabet
$\mathcal A$ represents the identity in the group given by
(\ref{ZP}) if and only if there exists a simply--connected planar
diagram $\Delta $ over (\ref{ZP}) such that $\phi (\partial \Delta
)\equiv W$ \cite[Ch. 5, Theorem 1.1]{LS}.

For every van Kampen diagram $\Delta $ over (\ref{ZP}) and any
fixed vertex $o$ of $\Delta $, there is a (unique) combinatorial
map $\gamma \colon Sk^{(1)} (\Delta )\to \Ga $ that preserves
labels and orientation of edges and maps $o$ to the vertex $1$ of
$\Ga $.

\vspace{3mm}

\noindent {\bf Hyperbolic spaces.} Here we briefly discuss some
properties of hyperbolic spaces used in this paper. For more
details we refer to \cite{BH,GH,Gro}.

One says that a metric space $M$ is {\it $\delta $--hyperbolic}
for some $\delta \ge 0$ (or simply {\it hyperbolic}) if for any
geodesic triangle $T$ in $M$, any side of $T$ belongs to the union
of the closed $\delta $--neighborhoods of the other two sides.

Recall that a path $p$ in a metric space is called {\it $(\lambda
, c)$--quasi--geodesic} for some $\lambda > 0$, $c\ge 0$,  if
$$dist(q_-, q_+)\ge \lambda l(q)-c$$ for any subpath $q$ of $p$.
Let $p$ be a path in a van Kampen diagram $\Delta $ over
(\ref{ZP}). We
need the following result about quasi--geodesics in hyperbolic
spaces (see, for example, \cite[Ch. III. H, Theorem 1.7]{BH}).

\begin{lem}\label{qg}
For any $\delta \ge 0$, $\lambda > 0$, $c\ge 0$, there exists a
constant $\kappa =\kappa (\delta , \lambda , c)$ with the
following property. If $M$ is a $\delta $--hyperbolic space and
$p, q$ are $(\lambda , c)$--quasi--geodesic paths in $M$ with same
endpoints, then $p$ and $q$ belong to the closed $\kappa
$--neighborhoods of each other.
\end{lem}

The next result can be easily derived from the definition of a
hyperbolic space by drawing the diagonal.

\begin{lem}\label{GeodQ}
Let $M$ be a $\delta $--hyperbolic metric space, $Q$ a geodesic
quadrangle in $M$. Then each side of $Q$ belongs to the closed
$2\delta $--neighborhood of the other three sides.
\end{lem}

From Lemma \ref{qg} and Lemma \ref{GeodQ}, we immediately obtain

\begin{cor}\label{qgq}
For any $\delta \ge 0$, $\lambda > 0$, $c\ge 0$, there exists a
constant $K=K(\delta , \lambda , c)$ with the following property.
Let $Q$ be a quadrangle in a $\delta$--hyperbolic space whose
sides are $(\lambda , c)$--quasi--geodesic. Then each side of $Q$
belongs to the closed $K$--neighborhood of the
union of the other three sides.
\end{cor}
\noindent Indeed it suffices to set $K=2(\kappa +\delta )$, where $\kappa $ is provided by Lemma \ref{qg}.

The proof of the next lemma is also straightforward (see \cite[Lemma 1.7]{Ols2}).

\begin{lem}\label{subseg}
Let $Q$ be a geodesic quadrangle with sides $a,b,c,d$ in a $\delta $-hyperbolic space. Suppose that $l(a)\ge 4\max \{ l(b),\, l(d)\} $. Then
there exist subsegments $p$ , $q$ of the sides $a$ and $b$, respectively, such that $\min \{ l(p), \, l(q)\} \ge \frac{7}{20} l(a) -8\delta $ and ${\rm dist\,} (p_{\pm }, q_\pm )\le 8\delta $.
\end{lem}

Passing from geodesic to quasi-geodesic quadrangles one can easily obtain the following. The proof is straightforward and consists of replacing the quasi-geometric quadrangle with a geodesic one having the same vertices, application of Lemma \ref{subseg}, and then Lemma \ref{qg}. We leave details to the reader.

\begin{cor}\label{subsegcor}
Let $Q$ be a $(\lambda , c)$--quasi-geodesic quadrangle with sides $a,b,c,d$ in a $\delta $-hyperbolic space, $\kappa =\kappa (\lambda , c)$ the constant provided by Lemma \ref{qg}. Suppose that $l(a)\ge (4\max \{ l(b),\, l(d)\} +c)/\lambda $. Then there exist subsegments $p$, $q$ of the sides $a$ and $b$, respectively, such that $\min \{ l(p), \, l(q)\} \ge \frac{7}{20} (\lambda l(a) -c) -8\delta - 2\kappa $ and ${\rm dist\,} (p_{\pm }, q_\pm )\le 8\delta + 2\kappa $.
\end{cor}

The following lemma is also well known (see, for example, \cite[CH.
III.H, Theorem 1.13]{BH}). Recall that a path in a metric space is
said to be $k$--local geodesic if any its subpath of length at
most $k$ is geodesic.

\begin{lem}\label{k-loc}
Let $r$ be a $k$--local geodesic in a $\delta $--hyperbolic metric
space for some $k>8\delta $. Then $r$ is $(1/3, 2\delta
)$--quasi--geodesic.
\end{lem}

The next lemma can be found in \cite[Lemma 25]{Ols1}.

\begin{lem}\label{N123}
There are positive constants $c_1=c_1(\delta )$ and $c_2=c_2(\delta )$ such that for any geodesic $r$-gon $P$ in a $\delta $-hyperbolic space, the following holds. Suppose that the set of sides of $P$ is divided into three subsets $N_1, N_2, N_3$ and $\sigma _i$ is the sum of lengths of sides from $N_i$, $i=1,2,3$. Assume that $\sigma _1>ar $ and $\sigma _3< 10^{-3}ar$ for some $a\ge c_2$.  Then there exist different sides $p_1\in N_1$ and $p_2\in N_1\cup N_2$ and subsegments $q_j$ of $p_j$, $j=1,2$, of lengths at least $10^{-3}a$ such that $$\max \{ {dist} ((q_1)_-, (q_2)_-),\, {dist} ((q_1)_+, (q_2)_+)\} \le c_1 .$$
\end{lem}

\vspace{3mm}

\noindent {\bf Relatively hyperbolic groups.} There are many
equivalent definitions of relatively hyperbolic groups (see
\cite{Bow,F,RHG} and references therein). In this paper we use the
isoperimetric characterization given in \cite{RHG}.

More precisely, let $G$ be a group, $\Hl $ a collection of
subgroups of $G$, $X$ a subset of $G$. We say that $X$ is a {\it
relative generating set of $G$ with respect to $\Hl $} if $G$ is
generated by $X$ together with the union of all $H_\lambda $. (In
what follows we $X$ to be symmetric.) In this situation the group
$G$ can be regarded as a quotient group of the free product
\begin{equation}
F=\left( \ast _{\lambda\in \Lambda } H_\lambda  \right) \ast F(X),
\label{F}
\end{equation}
where $F(X)$ is the free group with the basis $X$. Let $N$ denote
the kernel of the natural homomorphism $F\to G$. If $N$ is the
normal closure of a subset $\mathcal Q\subseteq N$ in the group
$F$, we say that $G$ has {\it relative presentation}
\begin{equation}\label{G}
\langle X,\; H_\lambda, \lambda\in \Lambda \; |\; \mathcal Q
\rangle .
\end{equation}
If $\sharp\, X<\infty $ and $\sharp\, \mathcal Q<\infty $, the
relative presentation (\ref{G}) is said to be {\it finite} and the
group $G$ is said to be {\it finitely presented relative to the
collection of subgroups $\Hl $.}

Set
\begin{equation}\label{H}
\mathcal H=\bigsqcup\limits_{\lambda\in \Lambda} (H_\lambda
\setminus \{ 1\} ) .
\end{equation}
Given a word $W$ in the alphabet $X\cup \mathcal H$ such that $W$
represents $1$ in $G$, there exists an expression
\begin{equation}
W=_F\prod\limits_{i=1}^k f_i^{-1}Q_i^{\pm 1}f_i \label{prod}
\end{equation}
with the equality in the group $F$, where $Q_i\in \mathcal Q$ and
$f_i\in F $ for $i=1, \ldots , k$. The smallest possible number
$k$ in a representation of the form (\ref{prod}) is called the
{\it relative area} of $W$ and is denoted by $Area^{rel}(W)$.

\begin{defn}
A group $G$ is {\it hyperbolic relative to a collection of
subgroups} $\Hl $ if $G$ is finitely presented relative to $\Hl $
and there is a constant $L>0$ such that for any word $W$ in $X\cup
\mathcal H$ representing the identity in $G$, we have $Area^{rel}
(W)\le L\| W\| $.
\end{defn}

In particular, $G$ is an ordinary {\it hyperbolic group} if $G$ is
hyperbolic relative to the trivial subgroup. An equivalent
definition says that $G$ is hyperbolic if it is generated by a
finite set $X$ and the Cayley graph $\Gamma (G, X)$ is hyperbolic.
In the relative case these approaches are not equivalent, but we
still have the following \cite[Theorem 1.7]{RHG}.

\begin{lem}\label{CG}
Suppose that $G$ is a group hyperbolic relative to a collection of
subgroups $\Hl $. Let $X$ be a finite relative generating set of
$G$ with respect to $\Hl $. Then the Cayley graph $\G $ of $G$
with respect to the generating set $X\cup \mathcal H$ is a
hyperbolic metric space.
\end{lem}

Observe also that the relative area of a word $W$ representing $1$
in $G$ can be defined geometrically via van Kampen diagrams. Let
$G$ be a group given by the relative presentation (\ref{G}) with
respect to a collection of subgroups $\Hl $.  We denote by
$\mathcal S$ the set of all words in the alphabet $\mathcal H$
representing the identity in the groups $F$ defined by (\ref{F}).
Then $G$ has the ordinary (non--relative) presentation
\begin{equation}\label{Gfull}
G=\langle X\cup\mathcal H\; |\;\mathcal S\cup \mathcal Q \rangle .
\end{equation}
A cell in van Kampen diagram $\Delta $ over (\ref{Gfull}) is
called a {\it $\mathcal Q$--cell} if its boundary is labeled by a
word from $\mathcal Q$. We denote by $N_\mathcal Q(\Delta )$ the
number of $\mathcal Q$--cells of $\Delta $. Obviously given a word
$W$ in $X\cup\mathcal H$ that represents $1$ in $G$, we have
$$
Area^{rel}(W)=\min\limits_{\phi (\partial \Delta ) \equiv W} \{
N_\mathcal Q (\Delta )\} ,
$$
where the minimum is taken over all van Kampen diagrams with
boundary label $W$.

\vspace{3mm}

\noindent {\bf $H_\lambda $--components.} Finally we are going to
recall an auxiliary terminology introduced in \cite{RHG}, which
plays an important role in our paper. Let $G$ be a group, $\Hl $ a
collection of subgroups of $G$, $X$ a finite generating set of $G$
with respect to $\Hl $, $q$ a path in the Cayley graph $\G $. A
subpath $p$ of $q$ is called an {\it $H_\lambda $--component} for
some $\lambda \in \Lambda $ (or simply a {\it component}) of $q$,
if the label of $p$ is a word in the alphabet $H_\lambda\setminus
\{ 1\} $ and $p$ is not contained in a bigger subpath of $q$ with
this property. Two $H_\lambda $--components $p_1, p_2$ of a path
$q$ in $\G $ are called {\it connected} if there exists a path $c$
in $\G $ that connects some vertex of $p_1$ to some vertex of
$p_2$ and ${\phi (c)}$ is a word consisting of letters from $
H_\lambda\setminus\{ 1\} $. In algebraic terms this means that all
vertices of $p_1$ and $p_2$ belong to the same coset $gH_\lambda $
for a certain $g\in G$. Note that we can always assume $c$ to have
length at most $1$, as every nontrivial element of $H_\lambda $ is
included in the set of generators.  An $H_\lambda $--component $p$
of a path $q$ is called {\it isolated } if no distinct $H_\lambda
$--component of $q$ is connected to $p$.

The lemma below is a simplification of from \cite[Lemma 2.27]{RHG}.

\begin{lem}\label{Omega}
Suppose that $G$ is a group that is hyperbolic relative to a
collection of subgroups $\Hl $, $X$ a finite generating set of $G$
with respect to $\Hl $. Then there exists a constant $K>0$ and a
finite subset $\Omega \subseteq G$ such that the following
condition holds. Let $q$ be a cycle in $\G $, $p_1, \ldots , p_k$
a set of isolated $H_\lambda $--components of $q$ for some
$\lambda\in \Lambda $, $g_1, \ldots , g_k$ the elements of $G$
represented by the labels of $p_1, \ldots , p_k$ respectively.
Then for any $i=1, \ldots , k$, $g_i$ belongs to the subgroup
$\langle \Omega \rangle \le G$ and the lengths of $g_i$ with
respect to $\Omega $ satisfy the inequality $$ \sum\limits_{i=1}^k
|g_i|_{\Omega }\le Kl(q).$$
\end{lem}

%%%%%%%%%%%%%%%%%%%%%%%%%%%%%%%%%%%%%%%%%%%%%%%%%%%%%%%%%%%%%%%%%

\section{Small cancellation conditions}

%%%%%%%%%%%%%%%%%%%%%%%%%%%%%%%%%%%%%%%%%%%%%%%%%%%%%%%%%%%%%%%%%

The main aim of this and the following four sections is to
generalize the small cancellation theory over hyperbolic groups
developed by Olshanskii in \cite{Ols2} to relatively hyperbolic
groups. The fact that the Cayley graph $\Gamma (G,X)$ of a
hyperbolic group $G$ generated by a finite set $X$ is a hyperbolic
metric space plays the key role in \cite{Ols2}. Lemma \ref{CG}
allows to extend this theory to the case of relatively hyperbolic
groups. However this extension is not straightforward as the
Cayley graph $\G $ defined in the previous section is not
necessary locally finite.

Roughly speaking, one can divide results about small cancellation conditions from \cite{Ols2}
into three classes. The first class consists of results about diagrams over presentations
satisfying small cancellation conditions, which do not use local finiteness of Cayley graphs at all. They can be stated and proved in our settings without any essential changes.  The main result of this kind is Lemma \ref{Gr0} stated below.

Proof of results from the second class do not use local finiteness of the Cayley graph either, but they do employ certain facts concerning geometric and algebraic properties of ordinary hyperbolic groups. These results can also be reproved with minor changes modulo the paper
\cite{RHG}, where the corresponding facts about relatively hyperbolic groups are contained.  For convenience of the reader we provide self-contained proofs in Sections 5 and 6.

Finally results of the third type explain how to chose words of some specific form satisfying small cancellation conditions. Unlike in ordinary small cancellation theory over a free group, verifying small cancellation conditions over hyperbolic groups is much harder and local finiteness of Cayley graphs is essentially used in \cite{Ols2}. Our approach here is different and is explained in Section \ref{wwsc}.

Given a set of words $\mathcal R$ in an alphabet $\mathcal A$, we
say that $\mathcal R$ is {\it symmetrized} if for any $R\in
\mathcal R$, $\mathcal R$ contains all cyclic shifts of $R^{\pm
1}$. Further let $G$ be a group generated by a set $\mathcal A$.
We say that a word $R$ is {\it $(\lambda , c)$--quasi--geodesic in
$G$}, if any path in the Cayley graph $\Gamma (G, \mathcal A)$
labeled $R$ is $(\lambda , c)$--quasi--geodesic.

\begin{defn}\label{piece}
Let $G$ be a group generated by a set $\mathcal A$, $\mathcal R$ a
symmetrized set of words in $\mathcal A$. For $\e
>0$, a subword $U$ of a word $R\in \mathcal R$ is called an {\it
$\e $--piece}  if there exists a word $R^\prime \in \mathcal R$
such that:

\begin{enumerate}
\item[(1)] $R\equiv UV$, $R^\prime \equiv U^\prime V^\prime $, for
some $V, U^\prime , V^\prime $; \item[(2)] $U^\prime = YUZ$ in $G$
for some words $Y,Z$ in $\mathcal A$ such that $\max \{ \| Y\| ,
\,\| Z\| \} \le \e $; \item[(3)] $YRY^{-1}\ne R^\prime $ in the
group $G$.
\end{enumerate}
Similarly, a subword $U$ of $R\in \mathcal R$ is called {\it an
$\e ^\prime $--piece}  if:
\begin{enumerate}
\item[($1^\prime $)] $R\equiv UVU^\prime V^\prime $ for some $V,
U^\prime , V^\prime $; \item[($2^\prime $)] $U^\prime =YU^{\pm
1}Z$ in the group $G$ for some $Y,Z$ satisfying $\max\{ \| Y\| ,
\| Z\| \}\le \e $.
\end{enumerate}
\end{defn}

Recall that a word $R$ in $\mathcal A$ is said to be {\it $(\lambda , c)$-quasi-geodesic}, if some (equivalently, any) path labeled by
$R$ in $\Gamma (G, \mathcal A)$ is $(\lambda , c)$-quasi-geodesic.

\begin{defn}\label{DefSC}
We say that the set $\mathcal R$ satisfies the {\it $C(\e , \mu ,
\lambda , c, \rho )$--condition} for some $\e \ge 0$, $\mu >0$,
$\lambda >0$, $c\ge 0$, $\rho >0$, if
\begin{enumerate}
\item[(1)] $\| R\| \ge \rho $ for any $R\in \mathcal R$;
\item[(2)] any word $R\in \mathcal R$ is $(\lambda ,
c)$--quasi--geodesic; \item[(3)] for any $\e $--piece of any word
$R\in \mathcal R$, the inequality $\max \{ \| U\| ,\, \| U^\prime
\| \} < \mu \| R\| $ holds (using the notation of Definition
\ref{piece}).
\end{enumerate}
Further the set $\mathcal R$ satisfies the {\it $C_1(\e , \mu ,
\lambda , c, \rho )$--condition} if in addition the condition
$(3)$ holds for any $\e ^\prime $--piece of any word $R\in
\mathcal R$.
\end{defn}

\begin{figure}
\hspace{3mm}\includegraphics{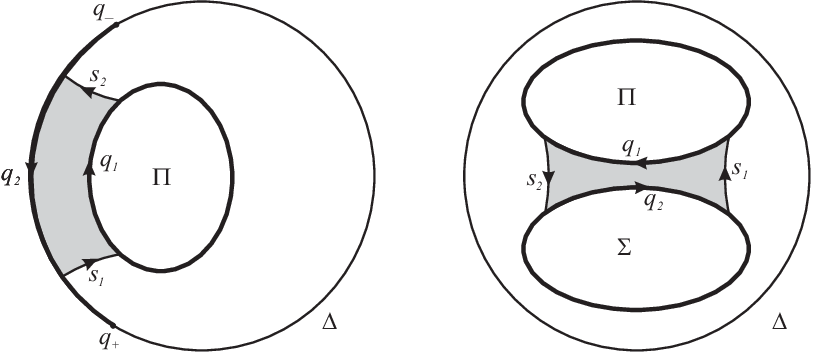}\\
\vspace{-3mm}
 \caption{Contiguity subdiagrams.}
\label{contig}
\end{figure}

Suppose that $G$ is a group defined by (\ref{ZP}). Given a set of
words $\mathcal R$, we consider the quotient group of $G$
represented by
\begin{equation}\label{quot}
G_1= \langle \mathcal A\; |\;
\mathcal O\cup \mathcal R\rangle .
\end{equation}
A cell in a van Kampen diagram over (\ref{quot}) is called an {\it
$\mathcal R$--cell} if its boundary label is a word from $\mathcal
R$.

Let $\Delta $ be a van Kampen diagram over (\ref{quot}), $q$ a subpath of $\partial \Delta $. Let also  $\Pi $, $\Sigma $ be
$\mathcal R$--cells of $\Delta $. Suppose that there is a simple
closed paths
\begin{equation}\label{pathp}
p=s_1q_1s_2q_2
\end{equation}
in $\Delta $, where $q_1$ is a subpath of
$\partial \Pi $, $q_2$ is a subpath of $q$ (or $\partial \Sigma $), and
\begin{equation}\label{sides}
\max \{l( s_1),\, l(s_2) \} \le \e
\end{equation}
for some constant $\e >0$. By $\Gamma $ we denote the subdiagram of
$\Delta $ bounded by $p$ (see Fig. \ref{contig}). If $\Gamma $
contains no $\mathcal R$--cells, we say that $\Gamma $ is an {\it
$\e $--contiguity subdiagram} (or simply a contiguity subdiagram if $\e $ is fixed) of $\Pi $ to the subpath $q$ of $\partial \Delta $ (or to $\Sigma $, respectively) and
$q_1$ is the {\it contiguity arc} of $\Pi $
to $q$ (respectively, $\Sigma $). The ratio $l(q_1)/l(\partial \Pi )$ is called the {\it contiguity degree}
of $\Pi $ to $q$ (respectively, $\Sigma $)
and is denoted by $(\Pi , \Gamma , q)$ (respectively, $(\Pi , \Gamma , \Sigma )$). In case $q=\partial \Delta $, we talk about contiguity subdiagrams, etc., of $\Pi $ to $\partial \Delta $.

A van Kampen diagram $\Delta $ over (\ref{quot}) is said to be
{\it reduced} if $\Delta $ has minimal number of $\mathcal
R$--cells among all diagrams over (\ref{quot}) having the same
boundary label.

The lemma below may be seen as a geometric reformulation of the small cancellation condition (see \cite[Lemma 5.2]{Ols2}). We provide the proof here for the sake of completeness.

\begin{lem}\label{O52}
Suppose that $\mathcal R$ satisfies the $C(\e,\mu, \lambda, c,\rho )$-condition for some values of the parameters. Let $\Delta $ be a reduced diagram over (\ref{quot}). Then for every $\e$-contiguity subdiagram $\Gamma $ of a cell $\Pi $ to another cell $\Sigma $, we have
\begin{equation}\label{geomSC}
\max \{ (\Pi , \Gamma , \Sigma ),\, (\Sigma , \Gamma , \Pi )\} <\mu .
\end{equation}
\end{lem}

\begin{proof}
Let $\partial \Gamma = s_1q_1s_2q_2$ as above (see Fig. \ref{contig}). Let $\partial \Pi=q_1a$ and $\partial \Sigma  =q_2b$. Then the word $U\equiv \phi (q_1)$ is an $\e$-piece. Indeed the first two conditions in Definition \ref{DefSC} are satisfied. If the third condition  is not, then $\phi (s_1) \phi (\partial \Pi )\phi (s_1)^{-1} =\phi (\partial \Sigma )$. That is, $\phi (s_1q_1as_1^{-1}bq_2)=1 $ in $G$. Hence we can cut the subdiagram of $\Delta $ bounded by $s_1q_1as_1^{-1}bq_2$ and fill the obtained hole with a diagram over (\ref{ZP}) without $\mathcal R$--cells reducing the number
of $\mathcal R$-cells by 2by $2$. This contradicts the assumption that $\Delta $ is reduced.  Thus $U$ is an $\e$-piece, and (\ref{geomSC}) follows from the $C(\e,\mu, \lambda, c,\rho )$--condition.
\end{proof}

When dealing with a diagram $\Delta $ over (\ref{quot}), it is convenient to consider the following transformations. Let $\Sigma $ be a subdiagram of $\Delta $ which contains no $\mathcal R$-cells, $\Sigma ^\prime $ another diagram over (\ref{ZP}) with $\phi (\partial \Sigma)=\phi (\partial \Sigma ^\prime )$. Then we can cut off $\Sigma $ and fill the obtained hole with $\Sigma ^\prime $. Note that this transformation does not affect $\phi(\partial \Delta )$ and the number of $\mathcal R$-cells in $\Delta $. If two diagrams over (\ref{quot}) can be obtained from each other by a sequence of such transformations, we call them {\it $\mathcal O$-equivalent}.

The following is an analogue of the well known Greendlinger Lemma. Recall that a path $p$ in $\Delta $ is called {\it $(\lambda , c)$-quasi-geodesic}, if its label is $(\lambda , c)$-quasi-geodesic, i.e., some (equivalently, any) path in $\Gamma (G, \mathcal A)$ with the same label is $(\lambda , c)$-quasi-geodesic.

\begin{lem}\label{Gr0}
Let $G$ be a group with a presentation $G=\langle \mathcal A\, |\, \mathcal O\rangle $. Suppose that the Cayley graph $\Gamma (G, \mathcal A)$ of $G$ is hyperbolic. Then for any $\lambda\in (0,1]$, $c\ge 0$, and $\mu \in (0, 1/16]$, there exist
$\e \ge 0$ and $\rho >0$ with the following property.
Let $\mathcal R$ be a symmetrized set of words in $\mathcal A$
satisfying the $C(\e, \mu , \lambda , c , \rho )$--condition,
$\Delta $ a reduced van Kampen diagram over the presentation
(\ref{quot}) whose boundary is $(\lambda , c)$--quasi--geodesic. Assume that $\Delta $ has at
least one $\mathcal R$--cell. Then there exists a diagram $\Delta ^\prime $ which is $\mathcal O$-equivalent to $\Delta $, an $\mathcal
R$--cell $\Pi $ in $\Delta ^\prime $, and an $\e$--contiguity subdiagram
$\Gamma $ of $\Pi $ to $\partial \Delta ^\prime $ such that
$$ (\Pi , \Gamma ,\partial \Delta ^\prime )> 1-13\mu .$$
\end{lem}

The proof of this lemma is actually given in \cite[Lemma 6.6]{Ols2} (see also the addendum in \cite{Ols3}) under the additional assumption that $\mathcal A$ is finite (i.e., the group $G$ is hyperbolic). However the finiteness of $\mathcal A$ was not used in the proof, so the same proof works in our case without any changes. Following the referee's recommendation, we provide a self-contained proof of Lemma \ref{Gr0} in Appendix.

%%%%%%%%%%%%%%%%%%%%%%%%%%%%%%%%%%%%%%%%%%%%%%%%%%%%%%%%%%%%%%%%%%

\section{Relative hyperbolicity of the quotient}

%%%%%%%%%%%%%%%%%%%%%%%%%%%%%%%%%%%%%%%%%%%%%%%%%%%%%%%%%%%%%%%%%%

Our next goal is to show that adding relations satisfying
sufficiently strong small cancellation conditions preserves
relative hyperbolicity. Throughout the rest of the paper (except the Appendix), let $G$ be a group hyperbolic
relative to a collection of subgroups $\Hl $, $X$ a finite
relative generating set of $G$ with respect to $\Hl $. We set
$\mathcal A=X\cup\mathcal H$ and $\mathcal O=\mathcal S\cup
\mathcal Q$, where $\mathcal S$ and $\mathcal Q$ are defined as in
(\ref{Gfull}). In the proof of the lemma below we follow the idea
from \cite[Lemma 6.7]{Ols2} with little changes.

\begin{lem}\label{G1}
For any $\lambda\in (0,1]$, $c\ge 0$, $N>0$, there exist $\mu
>0$, $\e \ge 0$, and $\rho >0$ such that for any finite
symmetrized set of words $\mathcal R$ satisfying the $C(\e , \mu ,
\lambda , c, \rho )$--condition, the following hold.
\begin{enumerate}

\item The group $G_1$ defined by (\ref{quot}) is hyperbolic
relative to the collection of images of subgroups $H_\lambda,
\lambda \in \Lambda $, under the natural homomorphism $G\to G_1$.
In particular, the Cayley graph of $G_1$ with respect to the
generating set $\mathcal A$ is hyperbolic.

\item The restriction of the natural homomorphism $\gamma \colon
G\to G_1$ to the subset of elements of length at most $N$ with
respect to the generating set $\mathcal A$ is injective.
\end{enumerate}
\end{lem}

\begin{proof}
Let us fix $\lambda ,c>0$ and sufficiently small $\mu< 1/16$. The exact value of $\mu $ required can easily extracted from the proof. If what follows, we assume that $\mu $ is small enough comparably to $\lambda $. Further we choose the constants $\e$ and $\rho $ according to Lemma \ref{Gr0}. Since the $\C $-condition becomes stronger as $\rho $ increases, we may increase $\rho $ if necessary without violating the conclusion of Lemma \ref{Gr0}. The exact lower bound can be easily extracted from our proof.

For a word $W$ in the alphabet $X\cup \mathcal H$ that represents
$1$ in $G_1$, we denote by $Area_1^{rel} (W)$ its relative area,
that is the minimal number $k$ in a representation of type
$$ W=_F\prod\limits_{i=1}^k f_iR_i^{\pm 1}f_i^{-1}$$ with the equality in the group
$F$ defined by (\ref{F}), where $R_i\in \mathcal R\cup \mathcal
Q$. Similarly, by $Area^{rel}(W) $ we denote the relative area of
a word $W$ representing $1$ in $G$, i.e., the minimal number $k$
in the above decomposition, where $R_i\in  \mathcal R$.

As $G$ is hyperbolic relative to $\Hl $, there exists a constant
$L>0$ such that for any word $W$ in $X\cup \mathcal H$
representing $1$ in $G$, we have $Area^{rel}(W)\le L\| W\| $.  To
prove the lemma it suffices to show that there is a constant
$\alpha >0$ such that for any word $W$ representing $1$ in $G_1$,
$Area_1^{rel}(W)\le \alpha \| W\| $. We are going to prove this
inequality for
$$\alpha =\frac{3L}{\lambda -27\mu }.$$

We proceed by induction on the length of $W$. Let $p$ be a path
labeled by $W$ in $\G $. Suppose first that $p$ is not $(1/2
,0)$--quasi--geodesic. Passing to a cyclic shift of $W$ if
necessary, we may assume that $p=p_0p_1$, where $p_0$ is a subpath
of $p$ such that $\dxh ((p_0)_-, (p_0)_+)< l(p_0)/2 $. Let $q$ be
a geodesic path in $\G $ that goes from $(p_0)_-$ to $(p_0)_+$.
Thus $l(q)< l(p_0)/2$. We denote by $U_0$, $U_1$, and $V$ the
labels of $p_0$, $p_1$, and $q$ respectively. Then
$$W=_F (U_0V^{-1})(VU_1) ,$$
where $U_0V^{-1}$ represents $1$ in $G$ and $VU_1$ represents $1$
in $G_1$. Obviously we have $$\| VU_1\| \le l(p)-l(p_0)+l(q)< \|
W\|-\frac{l(p_0)}{2} $$ and $$\| U_0V^{-1}\| \le l(p_0)+l(q) <
\frac{3l(p_0)}{2}.$$ Using the inductive hypothesis, we obtain
$$
\begin{array}{rl}
Area_1^{rel}(W) & \le Area_1^{rel}(VU_1)+Area^{rel}(U_0V^{-1})\\
& \\ & < \alpha \left( \| W\|-\frac12 l(p_0)\right) +\frac32 L
l(p_0)< \alpha \| W\|
\end{array}
$$
as $\alpha > 3L$.

Now suppose that $p$ is a $(1/2, 0)$--quasi--geodesic path in $\G
$. Since $C(\e , \mu , \lambda , c, \rho )$--condition implies
$C(\e , \mu , 1/2 , c, \rho )$ whenever $\lambda
>1/2$, it suffices to prove the lemma for $\lambda \le 1/2$.
So we may assume that $p$ is $(\lambda , c)$--quasi--geodesic as
well. Let $\Delta $ be a reduced diagram over (\ref{quot}) such
that $\partial \Delta $ is labeled $W$. Assume that $\Delta $ has
at least one $\mathcal R$--cell. (Otherwise the lemma is obviously
true.) Passing to an $\mathcal O$-equivalent diagram if necessary and using Lemma \ref{Gr0}, we can find an
$\mathcal R$--cell $\Pi $ and a contiguity subdiagram $\Gamma $ of
$\Pi $ to $\partial \Delta $ with $(\Pi , \Gamma , \partial \Delta
)>1-13\mu $. Let $\partial \Gamma =s_1q_1s_2q_2 $, where $q_1$
(respectively $q_2$) is a subpath of $\partial \Pi $ (respectively
$\partial \Delta $) and $\max \{ \| s_1\| ,\, \| s_2\| \} \le \e $. Let also $\partial \Pi =q_1u$ and
$\partial \Delta =wq_2$.

Note that perimeter of the subdiagram $\Xi $ of $\Delta $ bounded
by the path $s_2^{-1}us_1^{-1}w$ is smaller than perimeter of
$\Delta $ if $\rho $ is large enough and $\mu $ is close to zero.
Indeed as $\Gamma $ contains no $\mathcal R$--cells, we can regard
$s_1q_1s_2q_2$ as a cycle in the Cayley graph $\G $. Thus,
\begin{equation}\label{cont1}
\begin{array}{rl}
l(q_2)& \ge \dxh ((q_2)_-, (q_2)_+)\ge \dxh ((q_1)_-, (q_1)_+)-2\e
\ge \lambda l(q_1)-c-2\e \\ & \\ & \ge \lambda (1-13\mu )
l(\partial \Pi )-c-2\e \ge \lambda (1-13\mu )
\rho -c-2\e.
\end{array}
\end{equation}
(Recall that $l(\partial \Pi )\ge \rho $ by the $\C $-condition .) On the other hand,
\begin{equation}\label{cont2}
l(s_2^{-1}us_1^{-1})\le 2\e +13\mu l(\partial \Pi )\le 2\e +13\mu\rho.
\end{equation}
If $\rho $ is big enough and
$\mu $ is small enough, the right side of (\ref{cont1}) is
greater than the right side of (\ref{cont2}) and hence $l(\partial
\Xi )<l(\partial \Delta )$.

Therefore, by induction the total number $n_1$ of $\mathcal R$--
and $\mathcal Q$--cells in $\Xi $ is at most
\begin{equation}\label{n1}
\begin{array}{rl}
n_1& \le \alpha l(\partial \Xi )\le \alpha  \big( \| W\| -l(q_2) +
l (s_2^{-1}us_1^{-1})\big)
\\ & \\ & \le  \alpha \big( \| W\| -l(\partial \Pi )(\lambda
-26\mu ) +c +4\e \big) .
\end{array}
\end{equation}
Furthermore, as $q_2$ is $(1/2, 0)$--quasi--geodesic in $\G $, we
have
$$l(q_2)\le 2\dxh ((q_2)_-, (q_2)_+)\le 2(\dxh ((q_1)_-,
(q_1)_+)+2\e) \le 2l( \partial \Pi ) +4\e .$$ Therefore, the
perimeter of $\Gamma $ satisfies
$$
l(\partial \Gamma )\le 2\e +l(q_1)+ l(q_2)\le 3l(\partial \Pi )
+6\e .
$$
Hence we may assume that the number $n_2$ of $\mathcal Q$--cells
of $\Gamma $ is at most
\begin{equation}\label{n2}
n_2\le Ll(\partial \Gamma )\le 3L(l(\partial \Pi ) +2\e )\le
\alpha (l(\partial \Pi )(\lambda -27\mu ) +2\e ).
\end{equation}
Finally, combining (\ref{n1}) and (\ref{n2}), we obtain
$$
Area_1^{rel}(W)\le n_1+n_2\le \alpha ( \| W\| -\mu l(\partial \Pi
) +c+6\e )\le \alpha \| W\|
$$
whenever $\rho $ is big enough. This completes the
proof of relative hyperbolicity of $G_1$. The hyperbolicity of the
Cayley graph follows from Lemma \ref{CG}.

Note that if $\mu <1/30$, the inequality (\ref{cont1}) implies
$$
\| W\| \ge \frac12\lambda \rho -c-2\e
$$
for every non--trivial word $W$ which is geodesic in $G$ and
represents $1$ in $G_1$. Therefore, the second statement of the
lemma holds if we assume $\rho \ge 2(N+c+2\e )/\lambda $.
\end{proof}

%%%%%%%%%%%%%%%%%%%%%%%%%%%%%%%%%%%%%%%%%%%%%%%%%%%%%%%%%%%%%%%%%%

\section{Torsion in the quotient}

%%%%%%%%%%%%%%%%%%%%%%%%%%%%%%%%%%%%%%%%%%%%%%%%%%%%%%%%%%%%%%%%%%

We keep all assumptions and notation introduced in the beginning of the previous section.
Our next goal is to describe periodic elements in the quotient
group (\ref{quot}) of a relatively hyperbolic group $G$. To this
end we need an auxiliary result.

Recall that for an element
$g\in G$, the {\it
translation number} of $g$ with respect to $\mathcal A$ is defined
to be $$\tau _\mathcal A (g)=\lim\limits_{n\to \infty
}\frac{|g^n|_\mathcal A }{n}.$$ This limit always exists and is
equal to $\inf\limits_n (|g^n|_\mathcal A /n) $ \cite{GS}. The lemma
below can be found in \cite[Theorem 4.43]{RHG}.

\begin{lem}\label{cyc}
There exists $d>0$ such that for any hyperbolic element of
infinite order $g\in G$ we have $\tau _{X\cup \mathcal H} (g)\ge
d$.
\end{lem}

For our goals, even a stronger result is necessary.

\begin{lem}\label{Un}
There exist $1\ge \alpha >0$ and $a\ge 0$ with the following
property. Suppose that $g$ is a hyperbolic element of $G$ of
infinite order such that $g$ has the smallest length among all
elements of the conjugacy class $g^G$. Denote by $U$ a shortest
word in the alphabet $X\cup \mathcal H$ representing $g$. Then for
any $n\in \mathbb N$, any path in $\G $ labeled $U^n$ is
$(\alpha, a)$--quasi--geodesic.
\end{lem}

\begin{proof}
Recall that $\G $ is hyperbolic by Lemma \ref{CG}. First suppose
that $|g|_{X\cup \mathcal H}=\| U\| > 8\delta $, where $\delta $
is the hyperbolicity constant of $\G $. Since $g$ is a shortest
element in $g^G$ and $U$ is a shortest word representing $g$, there exists $k > 8\delta $ so
that the path $p$ labeled $U^n$ is a $k$--local geodesic in $\G$ for any $n$.  Therefore, by Lemma \ref{k-loc}, $p$
is $(1/3, 2\delta )$--quasi--geodesic.

Further if $|g|_{X\cup \mathcal H}=\| U\| \le 8\delta $, then for
any $n\in \mathbb N$, we have
$$ |g^n|_{X\cup \mathcal H}\ge n \inf\limits_i \left( \frac{1}{i} |g^i|_{X\cup \mathcal
H}\right) \ge n d\ge \frac{d}{8\delta }n|g|_{X\cup \mathcal H},$$
where $d$ is the constant provided by Lemma \ref{cyc}. Hence the
path $p$ labeled $U^n$ is $\left(\frac{d}{8\delta } , 8\delta
\right)$--quasi--geodesic. It remains to set $\alpha = \min\{
\frac13 ,\, \frac{d}{8\delta }\} $ and $a= 8\delta $.
\end{proof}

\begin{lem}\label{torsion}
For any $\lambda \in (0, 1]$, $c\ge 0$ there are $\mu >0$, $\e
\ge 0$, and $\rho >0$ such that the following condition holds. Suppose that $\mathcal R$ is a symmetrized set of words in $\mathcal A$
satisfying the $C_1 (\e, \mu , \lambda , c , \rho )$--condition.
Then every element of finite order in the group $G_1$ given by
(\ref{quot}) is the image of an element of finite order of $G$.
\end{lem}

\begin{proof}
Let us fix $\lambda , c>0$.
Observe that the $C_1 (\e, \mu,\lambda , c, \rho )$--condition becomes stronger as $\lambda $ increases and $c$ decreases. Hence it suffices to prove the lemma assuming that $\lambda < \alpha $ and $c>a$ for $\alpha $ and $a$ provided by Lemma \ref{Un}. Let us choose constants $1/16>\mu >0$, $\e >0$, $\rho >0$ such that
\begin{enumerate}
\item[($\ast$)] the conclusion of Lemma \ref{Gr0} holds.
\end{enumerate}
Again as in the proof of Lemma \ref{G1}, we note that the $C_1 (\e, \mu,\lambda , c, \rho )$--condition becomes stronger as $\mu $ decreases and $\e, \rho $ increase. Hence we may decrease $\mu$ and then increase $\rho $ and $\e$ if necessary without violating  ($\ast $). In particular, without loss of generality, we may  assume that
\begin{equation}\label{eps}
\e > 2\kappa + 8\delta ,
\end{equation}
where $\kappa =\kappa (\lambda, c)$ is the constant provided by Lemma \ref{qg} and $\delta $ is the hyperbolicity constant of $\G $. We fix $\e $ from now on.

Suppose that $g$ is an element of order $n>0$ in $G_1$ such that its preimage has infinite order in $G$. Assume also that $g$ has the smallest length among all elements from the conjugacy class $g^{G_1}$. Denote by $U$ a shortest word in the alphabet $X\cup \mathcal H$ representing $g$ in $G_1$. Then there exists a diagram $\Delta $ over (\ref{quot}) with boundary label $U^n$. By Lemma \ref{Un}, the label of $\partial \Delta $ is $(\lambda , c)$--quasi-geodesic (in $G$) and $\Delta $ contains at least one $\mathcal R$-cell. (For otherwise $\Delta $ is a diagram over (\ref{ZP}) and $U^n=1$ in $G$.) Note that passing from $\Delta $ to an $\mathcal O$-equivalent diagram does not affect $\partial \Delta $ and the number of $\mathcal R$-cells. Thus by Lemma \ref{Gr0} we can assume that there exists an $\mathcal R$--cell $\Pi $ in $\Delta $ and $\e $--contiguity subdiagram
$\Gamma $ of $\Pi $ to $\partial \Delta $ such that $(\Pi , \Gamma , \partial \Delta )> 1-13\mu .$

Let $\partial \Gamma =s_1q_1s_2q_2$ as on Fig. 1. Since $\Gamma $ has no $\mathcal R$-cells, we may think of $s_1q_1s_2q_2$ as a qadrangle in $\G $, where $q_i$ is $(\lambda , c)$-quasi-geodesic, and $l(s_i)\le \e $ for $i=1,2$. Therefore we have
\begin{equation}\label{lq2}
\begin{array}{rl}
l(q_2)&\ge \dxh ((q_2)_-, (q_2)_+) \ge \dxh ((q_1)_-, (q_1)_+)-l(s_1)-l(s_2) \\ &\\ &\ge \lambda l(q_1) -c -2\e \ge \lambda (1-13\mu) l(\partial \Pi ) -c-2\e \ge \lambda (1-13\mu) \rho -c-2\e .
\end{array}
\end{equation}
(Recall again that $l(\partial \Pi)\ge \rho$ by the $C_1 (\e, \mu,\lambda , c, \rho )$--condition.) In particular, choosing  $\rho $ large enough we can make  $l(q_2)$ as large as we want.

Passing from $U$ to a cyclic shift of $U^{\pm 1}$ if necessary, we may assume that $\phi (q_2)$ is a prefix of $U^n$ . We now have three cases to consider. Our goal is to show that neither of them is possible whenever $\mu $ is small enough and $\rho $ is large enough.

{\it Case 1.} Suppose that $l(q_2)\ge 4\| U\| /3$. This allows us to find two long disjoint equal subwords of $\phi (q_2)$. More precisely, we decompose  $\phi (q_2)$  as $\phi (q_2)\equiv WV_1WV_2$, where $$\lambda ^2 l(q_2)/5\le \| W\| \le \lambda ^2l(q_2)/4$$ and $$\| V_1\| >  l(q_2)/3.$$  Let $q_2=w_1v_1w_2v_2$ be the corresponding decomposition of the path $q_2$. Corollary \ref{qgq} applied to the quadrangle $s_1q_1s_2q_2$ implies that there is a point $o\in q_1$ such that $\dxh (o, (w_1)_+)\le K +\e$, where $K$ depends on $\lambda$, $c$, and $\delta $ only. Let $q_1^{-1}=r_1t$, where $(r_1)_+=o$. Thus we have
\begin{equation}\label{rw1}
\dxh ((r_1)_\pm , (w_1)_\pm )\le K+\e
\end{equation}
Similarly one can find a subpath $r_2^\xi $, $\xi =\pm 1$, of $q_1^{-1}$ such that
\begin{equation}\label{rw2}
\dxh ((r_2^\xi)_\pm , (w_2)_\pm )\le K+\e .
\end{equation}

We note that $r_1$ and $r_2$ are disjoint. Indeed otherwise $r_1$ passes through $(r_2)_-$ or $(r_2)_+$. For definiteness, assume that $(r_2)_-\in r_1$. Then we have
$$
\begin{array}{rl}
l(r_1) & \le \frac1\lambda (\dxh ((r_1)_-, (r_1)_+) +c)\le \frac1\lambda (l(w_1) +2\e + K+c)\\ &\\ & = \frac1\lambda (\|W\| +2\e + K+c)\le
\lambda  l(q_2)/4 + (2\e +K+c)/\lambda .
\end{array}
$$
On the other hand,
$$
\begin{array}{rl}
l(r_1)& \ge \dxh ((r_1)_-, (r_2)_-) \ge \lambda l(w_1v_1) -c-2\e -K\ge  \\ &\\ & \lambda \| V_1\| -c-2\e -K \ge \lambda l(q_2)/3 -c-2\e -K.
\end{array}
$$
These inequalities contradict each other if $l(q_2)$ is large enough. As explained before, the later condition can always be ensured by choosing sufficiently large $\rho $ (see (\ref{lq2})).

Thus we have a decomposition $q_1^{-1} =r_1t_1r_2^{\xi }t_2$, $\xi =\pm 1$. Let $\phi (q_1)^{-1} \equiv R_1T_1R_2T_2$ be the corresponding decomposition of the label.  By (\ref{rw1}) and (\ref{rw2}), we have $R_1=Y_1WZ_1$ and $R_2= Y_2W^{\pm 1} Z_2$ in $G$, where $\| Y_i\|, \| X_i\| \le K+\e $ for $i=1,2$. Hence $R_1=YR_2^{\pm 1 } Z$ in $G$, where  $\| Y\| , \| Z\| \le 2(K +\e)$. Without loss of generality, we may assume that words $Y$ and $Z$ are geodesic in $G$. Let $aybz$ be the corresponding rectangle in $\G$, where $\phi (a)\equiv R_1^{-1}$, $\phi (y)\equiv Y$, $\phi (b)\equiv R_2^{\pm 1}$, $\phi (z)\equiv Z$.  Recall that $R_1, R_2$ are subwords of the $(\lambda , c)$--quasi-geodesic word $\phi (\partial \Pi )$. Hence $a,b$ are $(\lambda , c)$--quasi-geodesic. Using (\ref{lq2}) and the triangle inequality, we obtain
\begin{equation}\label{lala}
\begin{array}{rl}
l (a)& =\| R_1\| \ge \dxh ((r_1)_-, (r_1)_-) \\ & \\ &\ge  \dxh ((w_1)_-, (w_1)_+) -l(s_1) -\dxh ((r_1)_+, (w_1)_+)\\ & \\ &  \ge \lambda l(w_1) -c-2\e -K= \lambda \| W\| -c - 2\e - K\\ & \\ &  \ge \lambda ^3 l(q_2) /5 -c-2\e -K \\ & \\ & \ge \lambda ^3[\lambda (1-13\mu) \rho -c-2\e ]/5 -c-2\e -K
\end{array}
\end{equation}
Clearly we can ensure $l(a)>(4\max\{ l(y), l(z)\} +c)/\lambda $ by taking large enough $\rho$. This allows us to apply Corollary \ref{subsegcor}. Thus we obtain subsegments $a^\prime , b^\prime $ of the sides $a$ and $b$, respectively, such that the distances between the corresponding ends of these subsegments is at most $8\delta +2\kappa $. Let $A=\phi (a^\prime )$, $B=\phi (b^\prime )$. Then $A=CB^{\pm 1}D$ in $G$, where $\| C\| , \| D\| \le 8\delta +2\kappa <\e$ by (\ref{eps}). Using this and (\ref{lala}), we obtain
$$
\min \{ \| A\|, \| B\| \}  \ge \frac{7}{20} (\lambda l(a) -c) -8\delta - 2\kappa \ge \mu\rho \ge \mu l(\partial \Pi )
$$
if $\mu $ is small enough and $\rho $ is large enough. Since $A$ and $B$ are disjoint subwords of $\phi (\partial \Pi )$, this contradicts the $C_1 (\e, \mu,\lambda , c, \rho )$--condition.

{\bf Case 2.} Suppose that $\| U\| \le l(q_2)\le 4\| U\| /3$, i.e. $\phi (q_2)\equiv UV$ for some (may be empty) word $V$, $\|V\|\le \| U\| /3$.
Note that $\phi (q_2)=\phi(s_2^{-1}us_1^{-1})$
in $G_1$, hence $g=\phi(s_2^{-1}us_1^{-1})V^{-1}$ in $G_1$. Since $U$ is the shortest word representing $g$ in $G_1$, we obtain
$$
\| U\| \le 2\e +13\mu l(\partial \Pi)  +\| V\| \le 2\e +13\mu l(\partial \Pi) + \| U\|/3.
$$
Consequently,
\begin{equation}\label{lule}
\| U \| \le 3(2\e +13\mu l(\partial \Pi))/2\le 3(2\e +13\mu \rho )/2.
\end{equation}
On the other hand, using (\ref{lq2}) we obtain
$$
\| U\| \ge 3l(q_2)/4 \ge 3(\lambda (1-13\mu) \rho -c-2\e )/4.
$$
The later inequality contradicts (\ref{lule}) whenever $\mu $ is small enough and $\rho $ is large enough.

{\bf Case 3.} Suppose that $ l(q_2)< \| U\| $, i.e., $\phi (q_2) $ is a subword of $U$. Again since $U$ is the shortest word representing $g$ in $G_1$, we have $\| \phi (q_2)\| \le \| Q\| $ for every word $Q$ such that $Q=\phi (q_2)$ in $G_1$. In particular, for $Q\equiv \phi (s_2^{-1}us_1^{-1})$, we obtain
\begin{equation}\label{Qge}
\| Q\| \ge  \| \phi (q_2)\| =l(q_2) \ge \lambda (1-13\mu) \rho -c-2\e .
\end{equation}
by (\ref{lq2}). On the other hand, we obviously have $\| Q\| \le 2\e +13\mu  l(\partial \Pi )<2\e +13\mu\rho$, which contradicts (\ref{Qge}) if $\mu $ is small enough and $\rho $ is large enough.
\end{proof}

%%%%%%%%%%%%%%%%%%%%%%%%%%%%%%%%%%%%%%%%%%%%%%%%%%%%%%%%%%%%%%%

\section{Words with small cancellations}\label{wwsc}

%%%%%%%%%%%%%%%%%%%%%%%%%%%%%%%%%%%%%%%%%%%%%%%%%%%%%%%%%%%%%%%

Recall that $G$ denotes a group hyperbolic
relative to a collection of subgroups $\Hl $, $X$ denotes a finite
relative generating set of $G$ with respect to $\Hl $. As above we
set $\mathcal A=X\cup \mathcal H.$ Our main goal here is to show
that a certain set of words over the alphabet $\mathcal A$
satisfies the small cancellation conditions described above.

More precisely, suppose that $W$ is a word satisfying the following
conditions.

\begin{enumerate}
\item[(W1)] $W\equiv xa_1b_1\ldots a_nb_n$ for some $n\ge 1$, where:

\item[(W2)] $x\in X\cup \{ 1\} $;

\item[(W3)] $a_1, \ldots , a_n$
(respectively $b_1, \ldots , b_n$) are elements of $H_\alpha $
(respectively $H_\beta $), where $H_\alpha \cap H_\beta =\{ 1\} $ ;

\item[(W4)] the elements $a_1, \ldots , a_n, b_1, \ldots , b_n$ do not belong to the set \begin{equation}\label{finiteset} \mathcal F = \mathcal F (\e ) =\{ g\in \langle \Omega \rangle\; : \; |g|_{\Omega} \le K(32\e +70) \} , \end{equation} where $\e $ is some non--negative integer and the set $\Omega $ and the constant $K$ are provided by Lemma \ref{Omega}.
\end{enumerate}

\begin{figure}
\unitlength=1mm \linethickness{0.4pt}
\begin{picture}(140,23)(-40,4)
\qbezier(65.23,7.95)(36.15,35.89)(6.01,7.95)

\put(5.83,7.95){\circle*{1}} \put(15.94,15.68){\circle*{1}}
\put(55.92,15.23){\circle*{1}} \put(50.28,18.65){\circle*{1}}
\put(22.03,18.95){\circle*{1}} \put(65.29,7.8){\circle*{1}}
\put(37.01,21.9){\vector(1,0){.07}}

\put(13,9){\makebox(0,0)[cc]{$p_1$}}
\put(60,9){\makebox(0,0)[cc]{$p_3$}}
\put(16.5,19.62){\makebox(0,0)[cc]{$s$}}
\put(54.5,19.3){\makebox(0,0)[cc]{$t$}}
\put(35,11.4){\makebox(0,0)[cc]{$e$}}
\put(35.64,25){\makebox(0,0)[cc]{$p_2$}}

\thicklines \linethickness{1pt}
\qbezier(15.76,15.61)(18.28,16.95)(22,18.88)
\qbezier(50.09,18.58)(53.74,16.65)(55.89,15.31)
\qbezier(22,18.88)(33.52,9.51)(50.09,18.58)
\put(54.07,16.6){\vector(2,-1){.07}}
\put(20,17.7){\vector(2,1){.07}}
\put(36.7,14.1){\vector(1,0){.07}}

\end{picture}

  \caption{}\label{ahm}
\end{figure}

Let $\mathcal {SW}$ denote the set of all subwords of cyclic
shifts of $W^{\pm 1}$. As in \cite{RHG}, we say that a path $p$ in
$\G $ is a {\it path without backtracking} if all components of
$p$ are isolated.

\begin{lem}\label{132}
Suppose $p$ is  a path in $\G $ such that $\phi (p)\in \mathcal
{SW}$. Then

1) $p$ is a path without backtracking.

2) $p$ is $(1/3, 2)$--quasi--geodesic.
\end{lem}

\begin{proof}
1) Suppose that $p=p_1sp_2tp_3$, where $s$ and $t$ are two
connected components. Passing to another pair of connected
components of $p$ if necessary, we may assume that $p_2$ is a path
without backtracking. For definiteness, we also assume that $s$
and $t$ are $H_\alpha $--components. Let $e$ denote a path of
length at most $1$ in $\G $ connecting $s_+$ to $t_-$ and labeled
by an element of $H_\alpha $ (see Fig. \ref{ahm}). It follows from
our choice of $W$ and the condition $H_\alpha \cap H_\beta =\{ 1\}
$ that $l(p_2)\ge 2$. Thus $p_2$ contains $m\ge l(p_2)/2\ge 1$
$H_\beta $--components, say $r_1, \ldots , r_m$, and all these
components are isolated components of the cycle $d=ep_2^{-1}$. Let
$g_1, \ldots , g_m $ be elements of $G$ represented by the labels
of $r_1, \ldots , r_m$. By Lemma \ref{Omega}, $g_i\in \langle
\Omega \rangle $, $i=1, \ldots , m$. According to (W4), we have
$$\sum\limits _{i=1}^m |g_i|_{\Omega } \ge 70Km\ge 35Kl(p_2)>
K(l(p_2)+1)= Kl(d).$$ This contradicts Lemma \ref{Omega}.

2) Since the set $\mathcal {SW}$ is closed under taking subwords,
it suffices to show that $\dxh (p_-, p_+)\ge  l(p)/3-2$. In case
$l(p)\le 6$ this is obvious, so we assume $l(p)> 6$. Suppose that
$\dxh (p_-, p_+)< l(p)/3-2$. Let $c$ denote a geodesic in $\G $
such that $c_-=p_-$ and $c_+=p_+$. Since $p$ is a path without
backtracking, any $H_\alpha $--component of $p$ is connected to at
most one $H_\alpha $--component of $c$. Obviously the path $p$
contains at least $l(p)/2-1$ $H_\alpha $--components. Therefore,
at least $$k=l(p)/2-1-l(c)>l(p)/2-1-(l(p)/3-2)> l(p)/6>1$$ of them
are isolated $H_\alpha $--components of the cycle $pc^{-1}$. Let
$f_1, \ldots , f_k$ be the elements of $G$ represented by these
components. Then as above we have $f_i\in \langle \Omega \rangle
$, $i=1, \ldots , k$. By (W4), we obtain
$$\sum\limits _{i=1}^k |f_i|_{\Omega } \ge 70Kk >
2Kl(p)>Kl(pc^{-1}).$$ This leads to a contradiction again.
\end{proof}

\begin{lem}\label{quad}
Suppose that $upv^{-1}q^{-1}$ is an an arbitrary quadrangle in $\G
$ satisfying the following conditions:

(a) $\phi (p)\equiv \phi (q^{-1})\in \mathcal {SW} $;

(b) $\max\{ l(u), \, l(v)\} \le \e $;

(c) $l(p)=l(q)\ge 6\e +22$.

\noindent Then the paths $p$ and $q$ have a common edge.
\end{lem}

The proof of Lemma \ref{quad} is based on the following two results. (We
keep the notation of Lemma \ref{quad} there.)

\begin{lem}\label{me1} Let us divide $p$ into three parts $p=p_1p_0p_2$ such that
\begin{equation}\label{p1p2}
l(p_1)=l(p_2)=3\e +6.
\end{equation}
Suppose that $s$ is a component of $p_0$. Then $s$ can not be
connected to a component of paths $u$ or $v$.
\end{lem}

\begin{proof}
Suppose that a component $s$ of $p_0$ is connected to a component
$t$ of $u$. Then $\dxh (s_+, t_+)\le 1$. Recall that the segment
$[p_-, s_+]$ of $p$ is $(1/3, 2)$--quasi--geodesic by Lemma
\ref{132}. Hence,
$$
\begin{array}{rl}
l(p_1) & < l([p_-, s_+])\le 3(\dxh (p_-, s_+) +2) \\ & \\
& \le 3(\dxh (p_-, t_+) +\dxh (t_+, s_+)+2) \\ &\\
& \le 3(l(u)-1+1+2)\le 3\e +6.
\end{array}
$$
However, this contradicts (\ref{p1p2}). Similarly $s$ can not be
connected to a component of $v$.
\end{proof}

\begin{lem}\label{me2}
Let $s_1, \ldots , s_k$ be consecutive components of $p_0$. Then
$q$ can be decomposed as $q=q_1t_1\ldots q_kt_kq_{k+1}$, where

1) $t_i$ is a component of $q$ connected to $s_i$, $i=1, \ldots ,
k$;

2) $q_i$ contains no components for $i=2,\ldots , k$, i.e. either
$\phi (q_i)\equiv x$ or $q_i$ is trivial.
\end{lem}

\begin{proof}
To prove the first assertion of the lemma we proceed by induction.
First let us show that $s_1$ is not isolated in
$d=upv^{-1}q^{-1}$.

Indeed assume $s_1$ is isolated in $d$. Suppose for definiteness
that $s_1$ is an $H_\alpha $--component. We consider the maximal
subpath $s$ of $p_0$ such that $s$ contains $s_1$ and all
$H_\alpha $--components of $s$ are isolated in $d$. By maximality
of $s$, either $s_-=(p_0)-$, or $s_-=r_+$ for a certain $H_\alpha
$--component $r$ of $p_0$ such that $r$ is not isolated in $d$.
(According to Lemma \ref{me1} this means that $r$ is connected to
an $H_\alpha $--component of $q$.) In the first case we denote by
$f_1$ the path $up_1$. In the second case, let $f_1$ be a path of
length $\le 1$ that connects an $H_\alpha $--component of $q$ to
$r$. In both cases we have $l(f_1)\le 4\e +6$. It follows from the
choice of $s$ that no $H_\alpha $--component of $s$ is connected
to an $H_\alpha $--component of $f_1$. Similarly we construct a
path $f_2$ such that $(f_2)_-\in q$, $(f_2)_+=s_+$, $l(f_2)\le 4\e
+6 $, and no $H_\alpha $--component of $s$ is connected to an
$H_\alpha $--component of $f_2$.

Clearly all $H_\alpha  $--components of $s$ are isolated in the
cycle $c=f_1sf_2^{-1}[(f_2)_-, (f_1)_-]$, where $[(f_2)_-,
(f_1)_-]$ is a segment  of $q^{\pm 1}$. We have
\begin{equation}\label{me21}
\dxh ((f_1)_-, (f_2)_-) \le l(f_1)+l(s)+l(f_2)\le 8\e + 12 +l(s).
\end{equation}
Consequently,
\begin{equation}\label{me22}
l([(f_1)_-, (f_2)_-])\le 3\dxh ((f_1)_-, (f_2)_-)+2)\le 24\e
+42+3l(s).
\end{equation}
Finally,
\begin{equation}\label{me23}
l(c)\le l(f_1)+l(s)+l(f_2)+l([(f_1)_-, (f_2)_-]) \le 32\e +54
+4l(s).
\end{equation}
Let $g_1, \ldots , g_m$ denote the elements represented by
$H_\alpha  $--components of $s$. Note that $l(s)\le 2m+2$.
Applying Lemma \ref{Omega}, we obtain $g_i\in \langle \Omega
 \rangle $, $i=1, \ldots , m$, and
\begin{equation}\label{me231}
\sum\limits_{i=1}^m |g_i|_{\Omega } \le Kl(c) \le K(32\e +54
+4l(s))\le K(32\e +62 +8m).
\end{equation}
Therefore, at least one of the elements $g_1, \ldots , g_m$ has
length at least
\begin{equation}\label{me232}
\frac{1}{m} K(32\e +62 +8m)\le K(32\e +70)
\end{equation}
that contradicts (W4). Thus $s_1$ is not isolated in $d$. By Lemma
\ref{me1} this means that $s_1$ is connected to an $H_\alpha
$--component $t_1$ of $q$.

Now assume that we have already found components $t_1, \ldots ,
t_i$ of $q$, $1\le i<k$, that are connected to $s_1, \ldots , s_i$
respectively. The inductive step is similar to the above
considerations. For definiteness, we assume that $s_i$ is an
$H_\alpha $--component. Then $s_{i+1}$ is an $H_\beta $--component
by the choice of $W$. We denote by $f_1$ a path of length $\le 1$
labeled by an element of $H_{\alpha }$ that connects $(t_i)_+$ to
$(s_i)_+$ (Fig. \ref{qfig}). If $s_{i+1}$ is isolated in the cycle
$c = f_1[(s_i)_+, p_+]v^{-1}[q_+, (t_i)_+]$, we denote by $s$ the
maximal initial subpath of the segment $[(s_i)_+, (p_0)_+]$ of
$p_0$ such that $s$ contains $s_{i+1}$ and all $H_\alpha
$--components of $s$ are isolated in $c$. As above, we can find a
path $f_2$ such that $(f_2)_-\in q$, $(f_2)_+=s_+$, $l(f_2)\le 4\e
+6 $, and no $H_\alpha  $--component of $s$ is connected to an
$H_\alpha $--component of $f_2$. The inequalities
(\ref{me21})--(\ref{me232}) remain valid and we arrive at a
contradiction in the same way. Thus $s_{i+1}$ is not isolated in
$c$, i.e., it is connected to a component $t_{i+1}$ of the segment
$[(t_i)_+, q_+]$ of $q$. This completes the inductive step.

\begin{figure}

\unitlength 1mm % = 2.85pt
\linethickness{0.4pt}
\begin{picture}(124.75,36)(-13,50)
\qbezier[1000](5.13,81.67)(60.55,73.89)(124.1,81.67)
\qbezier[1000](5.13,54.09)(60.55,61.87)(124.1,54.09)

\put(7.2,68.38){\vector(0,1){.07}}\qbezier(5.13,54)(9.13,68.94)(5.13,81.63)

\put(121.33,68.38){\vector(0,1){.07}}\qbezier(124.13,54.13)(118.38,68.88)(124.13,81.63)

\thicklines

\put(42,78.4){\vector(1,0){.07}}\qbezier(33.25,78.75)(42.94,78.25)(48.88,78)

\put(33.25,78.63){\circle*{1}} \put(60,78){\circle*{1}}
\put(77,78){\circle*{1}} \put(90.13,78.75){\circle*{1}}
\put(90.13,57.25){\circle*{1}} \put(30.13,56.75){\circle*{1}}
\put(19.75,55.88){\circle*{1}} \put(48.88,77.88){\circle*{1}}
\put(26,56.38){\vector(1,0){.07}}\qbezier(19.63,55.88)(27.13,56.5)(30.13,56.8)

\put(68.5,77.94){\vector(1,0){.07}}\multiput(59.88,77.88)(4.3125,.0313){4}{\line(1,0){4.3125}}
\put(42.1,67.5){\vector(1,1){.07}}\qbezier(30,56.75)(44.38,67.69)(48.75,77.88)

\thinlines

\put(88.1,68.38){\vector(0,1){.07}}\qbezier(90,57.25)(86,68.69)(90,78.88)

\put(4,68.13){\makebox(0,0)[cc]{$u$}}
\put(124.5,69){\makebox(0,0)[cc]{$v$}}
\put(11.88,52.38){\makebox(0,0)[cc]{$q$}}
\put(12.75,83.13){\makebox(0,0)[cc]{$p$}}
\put(25.88,53){\makebox(0,0)[cc]{$t_i$}}
\put(42,81){\makebox(0,0)[cc]{$s_i$}}
\put(49.38,80.2){\makebox(0,0)[cc]{$s_-$}}
\put(90,81.3){\makebox(0,0)[cc]{$s_+$}}
\put(67.75,80.8){\makebox(0,0)[cc]{$s_{i+1}$}}
\put(44.4,66){\makebox(0,0)[cc]{$f_1$}}
\put(85.25,67.63){\makebox(0,0)[cc]{$f_2$}}
\put(65.13,67.13){\makebox(0,0)[cc]{$c$}}
\end{picture}

\unitlength 1mm % = 2.85pt
\linethickness{0.4pt}
\ifx\plotpoint\undefined\newsavebox{\plotpoint}\fi % GNUPLOT compatibility
  \caption{}\label{qfig}
\end{figure}

Let us prove the second assertion of the lemma. Suppose that $s_i,
s_{i+1}$ are two subsequent components of $p_0$ and such that
$q_{i+1}$ contains at least one component (say, an $H_\alpha
$--component). As above for definiteness we assume that $s_i$
(respectively $s_{i+1}$) is an $H_\alpha $--component
(respectively $H_\beta $--component). Let $e_1$, $e_2$ be the
paths of lengths $\le 1$ labeled by elements of $H_\alpha $ and
$H_\beta $ respectively such that $(e_1)_+=(s_i)_+$,
$(e_1)_-=(q_{i+1})_-$, $(e_2)_+=(s_{i+1})_-$,
$(e_2)_-=(q_{i+1})_+$ (see Fig. \ref{3}). As $q$ is a path without
backtracking, each $H_\alpha $--component of $q_{i+1}$ is isolated
in the cycle $e= q_{i+1} e_2[(s_{i+1})_-, (s_{i})_+]e_1^{-1}$,
where $[(s_{i+1})_-, (s_{i})_+]$ is a segment of $p^{-1}$. Notice
that $l([(s_{i+1})_-, (s_{i})_+])\le 1$ as $s_i$ and $s_{i+1}$ are
subsequent components of $p_0$. We denote by $f_1, \ldots , f_m$
the elements represented by $H_\alpha $--components of $q_{i+1}$.
By Lemma \ref{Omega}, we have $f_j\in \langle \Omega \rangle $,
$j=1, \ldots , m$, and
$$
\sum\limits_{j=1}^m |f_j|_{\Omega } \le Kl(e) \le
K(3+l(q_{i+1}))\le K(3+2m+2)\le 7mK.
$$
This contradicts (W4) again.
\end{proof}

\begin{figure}
\unitlength 1mm
  \begin{picture}(124.75,36)(-15,50)
\qbezier[1000](5.13,81.67)(60.55,73.89)(124.1,81.67)
\qbezier[1000](5.13,54.09)(60.55,61.87)(124.1,54.09)

\put(7.13,68.38){\vector(0,1){.07}}\qbezier(5.13,54)(9.13,68.94)(5.13,81.63)

\put(121.33,68.38){\vector(0,1){.07}}\qbezier(124.13,54.13)(118.38,68.88)(124.13,81.63)

\thicklines

\put(42,78.4){\vector(1,0){.07}}\qbezier(33.25,78.75)(42.94,78.25)(48.88,78)

\put(33.25,78.63){\circle*{1}} \put(60,78){\circle*{1}}
\put(77,78){\circle*{1}} \put(80.25,57.63){\circle*{1}}
\put(90,57.13){\circle*{1}} \put(30.13,56.75){\circle*{1}}
\put(19.75,55.88){\circle*{1}} \put(48.88,77.88){\circle*{1}}

\put(26,56.38){\vector(1,0){.07}}\qbezier(19.63,55.88)(27.13,56.5)(30.13,56.8)

\put(68.5,77.94){\vector(1,0){.07}}\multiput(59.88,77.88)(4.3125,.0313){4}{\line(1,0){4.3125}}

\put(4,68.13){\makebox(0,0)[cc]{$u$}}
\put(124.5,69){\makebox(0,0)[cc]{$v$}}
\put(11.88,52.38){\makebox(0,0)[cc]{$q$}}
\put(12.75,83.13){\makebox(0,0)[cc]{$p$}}
\put(25.88,53){\makebox(0,0)[cc]{$t_i$}}
\put(42,81){\makebox(0,0)[cc]{$s_i$}}

\put(67.75,80.63){\makebox(0,0)[cc]{$s_{i+1}$}}

\put(85.88,57.38){\vector(1,0){.07}}\qbezier(80.13,57.63)(86.75,57.31)(89.88,57.25)

\put(58.38,57.88){\vector(1,0){.07}}\put(58,57.88){\line(1,0){.375}}

\put(42.38,67.13){\vector(1,1){.07}}\qbezier(30.13,56.63)(45.19,67)(49,77.88)

\put(65.13,67.38){\vector(-1,1){.07}}\qbezier(80.13,57.63)(60.25,66.94)(59.88,78)

\put(84.5,53.88){\makebox(0,0)[cc]{$t_{i+1}$}}
\put(57.13,54.5){\makebox(0,0)[cc]{$q_{i+1}$}}
\put(45,65.5){\makebox(0,0)[cc]{$e_1$}}
\put(63.75,65.5){\makebox(0,0)[cc]{$e_2$}}
\put(54.5,68.88){\makebox(0,0)[cc]{$e$}}
\end{picture}
  \caption{}\label{3}
\end{figure}

\begin{proof}[Proof of Lemma \ref{quad}]
We keep the notation introduced in Lemma 5.3 and Lemma 5.4. Let
also $$p_0=p_1s_1\ldots p_ks_kp_{k+1}$$ for some (may be trivial)
subpaths $p_1, \ldots p_{k+1}$ of $p_0$.

According to (\ref{p1p2}) and condition c) of Lemma \ref{quad}, we
have $$ l(p_0)\ge 6\e +22 -l(p_1)-l(p_2)\ge 10.$$ Since $s_1,
\ldots , s_k$ are subsequent components and $\phi (p_0)$ is a
subword of a cyclic shift of $W^{\pm 1}$, at most one of the paths
$p_2, \ldots , p_{k}$ is non--trivial. Therefore, there are at
least $5$ subsequent components $s_{i}, \ldots , s_{i+4}$, such
that $p_{i+1}, \ldots , p_{i+4}$ are trivial. Without loss of
generality we may assume $i=1$. Similarly, by Lemma \ref{me2}, we
can find at least three subsequent components among $t_1, \ldots
t_5$, say $t_1, t_2, t_3$, such that $(t_1)_+=(t_2)_-$ and
$(t_2)_+=(t_3)_-$. Let $w$ be an element represented by the label
of any path that goes from $(t_1)_+=(t_2)_-$ to $(s_1)_+=(s_2)_-$.
For definiteness we assume that $t_1$ and $s_1$ are $H_\alpha
$--components. Since $t_1$ and $s_1$ are connected, we have $w\in
H_\alpha $. On the other hand, the $H_\beta $--components $t_2$
and $s_2$ are also connected. Hence $w\in H_\beta $. Thus $w\in
H_\alpha \cap H_\beta =\{ 1\} $, i.e. the vertices $(t_2)_-$ and
$(s_2)_-$ coincide. Similarly the vertices $(t_2)_+$ and $(s_2)_+$
coincide. In particular, $t_2$ and $s_2$ are edges labeled by the
same element of $H_\beta $, i.e., $t_2$ and $s_2$ coincide.
\end{proof}

Now we are ready to prove the main result of this section.

\begin{thm}\label{w}
Suppose that $W$ is a word in $\mathcal A$ satisfying the
conditions (W1)--(W4) and, in addition, $a_i\ne a_j^{\pm 1}$,
$b_i\ne b_j^{\pm 1}$ whenever $i\ne j$ and $a_i\ne a_i^{-1}$,
$b_i\ne b_i^{-1}$, $i,j\in \{ 1, \ldots , n\} $. Then the set
$\mathcal W$ of all cyclic shifts of $W^{\pm 1} $ satisfies the
$C_1(\e , \frac{3\e +11}{n} ,\frac13, 2, 2n+1)$ small cancellation
condition.
\end{thm}

\begin{proof}
The first two conditions from Definition \ref{DefSC} follow from
the choice of $W$ and Lemma \ref{132}. Suppose that $U$ is an $\e
$--piece of a word $R\in \mathcal {W}$. Assume that $\max \{ \|
U\| , \| U^\prime \| \} \ge \mu \| R\| $ for $\mu = \frac{3\e
+11}{n} $, that is, $$\max \{ \| U\| , \| U^\prime \| \} \ge
\frac{3\e +11}{n} (2n+1) > 6\e +22.$$ (Here and below we use the
notation of Definitions \ref{DefSC} and \ref{piece}.) Without loss
of generality we may assume that $\| U\| \ge 6\e +22$. By the
definition of an $\e $--piece, there is a quadrangle
$upv^{-1}q^{-1}$ in $\G $ satisfying conditions (a)--(c) of Lemma
\ref{quad} and such that labels of $p$ and $q$ are $U$ and
$U^\prime $ respectively.

Let $e$ be the common edge of $p$ and $q$. Then we have $$R\equiv
U_1\phi (e)U_2V$$ and $$R^\prime \equiv U_1^\prime\phi
(e)U_2^\prime V^\prime ,$$ where $U_1\phi (e)U_2\equiv U$ and
$U^\prime \equiv U_1^\prime\phi (e)U_2^\prime $. Since $\phi (e)$
appears in $W^{\pm 1}$ only once, $R$, $R^\prime $ are cyclic
shifts of the same word $W^{\pm 1}$ and
$$U_2VU_1\equiv U_2^\prime V^\prime U_1^\prime .$$ Note also that
$$Y=U_1^\prime U_1^{-1} $$ in $G$ as $Y^{-1}U_1^\prime U_1^{-1} $ is
a label of a cycle in $\G $. Therefore, the following equalities
hold in the group G:
$$
YRY^{-1}= U_1^\prime U_1^{-1} U_1\phi (e)U_2V U_1(U_1^\prime
)^{-1} = U_1^\prime \phi (e) U_2^\prime V^\prime U_1 ^\prime
(U_1^\prime )^{-1} =R^\prime
$$
that contradicts the third condition from Definition \ref{DefSC}.

Similarly, if $U$ is an $\e ^\prime $--piece, then $R\equiv
UVU^\prime V^\prime $ for some $U, U^\prime ,V, V^\prime $, where
both subwords $U$ and $U^\prime $ contain a certain letter from
$X\cup\mathcal H$. In this case we arrive at a contradiction again
as any letter $a\in X\cup\mathcal H$ appears in $R$ only once, and
if $a$ appears in $R$, then $a^{-1}$ does not.
\end{proof}

Finally we note that the condition $x\in X\cup \{ 1\} $ in Theorem
\ref{w} is not really restrictive since we can always add any element $x\in G$ to the
set $X$ without violating  relative hyperbolicity.

%%%%%%%%%%%%%%%%%%%%%%%%%%%%%%%%%%%%%%%%%%%%%%%%%%%%%%%%%%%%%%%%%%

\section{Suitable subgroups and quotients}

%%%%%%%%%%%%%%%%%%%%%%%%%%%%%%%%%%%%%%%%%%%%%%%%%%%%%%%%%%%%%%%%%%

Throughout this section, we keep the assumption that $G$ is
hyperbolic relative to a collection of subgroups $\Hl $. The proof
of Lemma 2.3 is based on  the following auxiliary result.

\begin{lem}\label{f1f2}
Suppose that for some
$\lambda , \mu \in \Lambda $, $\lambda \ne \mu $, $H_\lambda $ and $H_\mu $ contain elements of infinite order and $H_\lambda \cap H_\mu=\{ 1\} $. Then there are $f\in H_\lambda $, $g\in H_\mu $ such that $fg$ is a hyperbolic element of infinite order and $E_G(fg)=\langle fg \rangle $.
\end{lem}

\begin{proof}
We set $\e = 2(\kappa + \delta )$, where $\kappa =\kappa (\delta
,1/3 , 2 )$ is the constant provided by Lemma \ref{qg} and $\delta
$ is the hyperbolicity constant of $\G $. It is convenient to assume that $\kappa , \delta \in \mathbb N$.

Let $\mathcal F=\mathcal
F (\e )$ be the set defined by (\ref{finiteset}). Since $\Omega $
is finite, $\mathcal F$ is finite, and hence there are elements $f\in
H_\lambda \setminus \mathcal F$ and $g\in H_\mu \setminus \mathcal
F$ of infinite order. In particular,
\begin{equation}\label{invol}
f^2\ne 1,\;\;\; g^2\ne 1.
\end{equation}
Note that the word $W=(fg)^m$ satisfies conditions (W1)--(W4).
(Here $f$ and $g$ are regarded as letters of $\mathcal H$.) By
Lemma \ref{132}, for any $m\in \mathbb N$, the word $(fg)^m$ is
$(1/3, 2)$--quasi--geodesic in $\G $. In particular, $fg$ is a
hyperbolic element. Indeed otherwise the length $|(fg)^m|_{X\cup
\mathcal H}$ would be bounded uniformly on $m$.

Now suppose that $a\in E_G(fg)$. Then by Theorem \ref{E(g)},
$a(fg)^ma^{-1}=(fg)^{\pm m }$ for some $m\in \mathbb N$. Passing
to a multiple of $m$ if necessary, we may assume that
\begin{equation} \label{m}
m\ge 3|a|_{X\cup \mathcal H} +12(\kappa +\delta) +20.
\end{equation}

Let $a_1b_1a_2b_2$ be a quadrangle in $\G $ such that $a_1$, $a_2$
are geodesic, the labels $\phi (a_1)=\phi (a_2^{-1})$ represent
$a$ in $G$, and $\phi (b_1)=\phi (b_2^{-1})\equiv (fg)^{\pm m}$.
Let also $b_1=b_1^\prime pb_1^{\prime\prime } $ (Fig.
\ref{E(fg)}), where
\begin{equation}\label{b1}
l(b_1^\prime )=l(b_1^{\prime\prime }) =3(l(a_1)+2(\kappa +\delta )
+3 ).
\end{equation}
As $b_1$ is $(1/3, 2)$--quasi--geodesic, we have
\begin{equation}\label{ss1}
\dxh (p_\pm , (b_1)_\pm ) \ge \frac13 l(b_1^\prime ) -2 \ge l(a_1)
+2(\kappa +\delta )+1.
\end{equation}
Further by Corollary \ref{qgq}, there is a point $s\in a_1\cup
b_2\cup a_2 $ such that $\dxh (s, p_-)\le 2(\kappa +\delta )$. If
$s\in a_1$, then we have $$\dxh (p_-, b_-)\le \dxh (p_-, s)+\dxh
(s, b_-)\le 2(\kappa +\delta )+l(a_1)$$ that contradicts
(\ref{ss1}). For the same reason $s$ can not belong to $a_2$. Thus
$s\in b_2$. Without loss of generality, we may assume that $s$ is
a vertex of $\G$. Similarly there exists a vertex $t\in b_2$ such
that $\dxh (t, p_+)\le 2(\kappa +\delta )$. Let $u,v$ be geodesics
in $\G $ connecting $s$ to $p_-$ and $t$ to $p_+$ respectively and
let $q$ denote the segment $[s,t]$ of $b_2^{-1}$.

According to (\ref{m}) and (\ref{b1}), we have
$$
l(p)\ge 2m - 6(l(a_1)+2(\kappa +\delta )+3)\ge 12 (\kappa +\delta
)+22.
$$
Hence we may apply Lemma \ref{quad} for the quadrangle
$upv^{-1}q^{-1}$ and $\e =2(\kappa +\delta )$. Thus there exists a
common edge $e$ of $p$ and $q$. In particular, this and (\ref{invol}) imply that
$a(fg)^ma^{-1}=(fg)^{m }$ (not $(fg)^{-m}$).

There are two possibilities for labels of the segments $[(b_1)_-,
e_-]$ and $[(b_2)_+, e_-]$ of $b_1$ and $b_2^{-1}$ respectively.
Namely both these labels are either of the form $(fg)^n$ (possibly
for different $n$) or of the form $(fg)^kf$. In both cases
$a=(fg)^l$ for a certain $l$ as labels of $a_1$ and $[(b_2)_+,
e_-][(b_1)_-, e_-]^{-1}$ represent the same element of $G$. Thus
$a\in \langle fg \rangle $ and $E_G(fg)=\langle fg\rangle $.
\end{proof}

\begin{figure}
\unitlength 1mm % = 2.85pt
\linethickness{0.4pt}
\ifx\plotpoint\undefined\newsavebox{\plotpoint}\fi % GNUPLOT compatibility
\begin{picture}(116.48,25)(-18,35)
\put(55.95,44.9){\vector(1,0){.07}}\put(45.96,44.9){\line(1,0){19.976}}
\qbezier(65.94,44.9)(100.67,49.06)(115.97,57.81)
\qbezier(67.04,44.89)(96.1,43.18)(115.95,35.82)
\qbezier(45.93,44.74)(21.63,44.37)(2.97,35.97)
\qbezier(45.93,44.74)(18.88,47.05)(3.12,57.97)
\put(45.78,44.89){\circle*{1}} \put(66.9,44.89){\circle*{1}}
\put(97.96,50.54){\circle*{1}} \put(97.96,41.03){\circle*{1}}
\put(115.8,35.82){\circle*{1}} \put(115.95,57.68){\circle*{1}}
\put(3.12,35.97){\circle*{1}} \put(3.57,57.68){\circle*{1}}
\put(18.14,50.39){\circle*{1}} \put(17.84,41.03){\circle*{1}}

\put(18.73,46.53){\vector(0,1){.07}}\qbezier(17.84,41.03)(19.47,47.35)(18.14,50.39)
\put(97.5,46.38){\vector(0,1){.07}}\qbezier(97.96,41.03)(96.85,46.97)(97.81,50.54)
\put(4.91,46.38){\vector(0,1){.07}}\qbezier(3.12,35.82)(6.61,46.01)(3.27,57.68)
\put(10.7,53.51){\vector(2,-1){.07}}\multiput(10.41,53.66)(.05946,-.02973){5}{\line(1,0){.05946}}
\put(105.1,52.94){\vector(3,1){.07}}

\put(12,55.89){\makebox(0,0)[cc]{$b_1^\prime $}}
\put(105.5,55.5){\makebox(0,0)[cc]{$b_1^{\prime\prime }$}}
\put(96.3,53){\makebox(0,0)[cc]{$p_+$}}
\put(21,53){\makebox(0,0)[cc]{$p_-$}}
\put(55,47.42){\makebox(0,0)[cc]{$e$}}
\put(22,38.2){\makebox(0,0)[cc]{$s=q_-$}}
\put(95,38.5){\makebox(0,0)[cc]{$t=q_+$}}
\put(2.23,45.64){\makebox(0,0)[cc]{$a_1$}}

\put(113.8,45.04){\vector(0,-1){.07}}\qbezier(115.8,57.83)(111.64,43.18)(115.8,35.97)

\put(117,46){\makebox(0,0)[cc]{$a_2$}}
\put(100,45.78){\makebox(0,0)[cc]{$v$}}
\put(16.2,46.08){\makebox(0,0)[cc]{$u$}}
\end{picture}

  \caption{}\label{E(fg)}
\end{figure}

\begin{proof}[Proof of Lemma \ref{non-com}]
Let $f_1, f_2\in H^0$ be non--commensurable elements of $H$ such
that $E_G(f_1)\cap E_G(f_2)=\{ 1\} $. By Theorem \ref{E(g)}, $G$
is hyperbolic relative to the collection $\Hl\cup E_G(f_1)\cup
E_G(f_2)$.

We construct a sequence of desired elements $h_1, h_2, \ldots $ by
induction. By Lemma \ref{f1f2}, there are $f\in E_G(f_1)$, $g\in
E_G(f_2)$ such that the element $h_1=fg$ is hyperbolic ( with
respect to the collection $\Hl\cup E_G(f_1)\cup E_G(f_2)$) and
$E_G(h_1)=\langle h_1\rangle $. Theorem \ref{E(g)} implies that
$G$ is hyperbolic relative to the collection $\Hl\cup E_G(f_1)\cup
E_G(f_2)\cup E_G(h_1)$. Further we construct a hyperbolic (with
respect to $\Hl\cup E_G(f_1)\cup E_G(f_2)\cup E_G(h_1)$) element
$h_2$ as a product of an element of $E_G(f_1)$ and an element of
$E_G(f_2)$ as above. As $h_2$ is hyperbolic, it is not
commensurable with $h_1$. Applying Theorem \ref{E(g)} again we
join $E_G(h_2)$ to the collection of subgroups with respect to
which $G$ is hyperbolic and so on. Continuing this procedure, we
get what we need.
\end{proof}

To prove Theorem \ref{glue}, we need the following two
observations. The first one is a particular case of Theorem 2.40
from \cite{RHG}.

\begin{lem}\label{exhyp}
Suppose that a group $G$ is hyperbolic relative to a collection of
subgroups $\Hl \cup \{ S_1, \ldots , S_m\} $, where $S_1, \ldots ,
S_m $ are finitely generated and hyperbolic in the ordinary
(non--relative) sense. Then $G$ is hyperbolic relative to $\Hl $.
\end{lem}

The next lemma is a particular case of Theorem 1.4 from \cite{RHG}

\begin{lem}\label{malnorm}
Suppose that a group $G$ is hyperbolic relative to a collection of
subgroups $\Hl $. Then

\begin{enumerate}
\item[(a)] For any $g\in G$ and any $\lambda , \mu \in \Lambda $,
$\lambda \ne \mu $, the intersection $H_\lambda ^g\cap H_\mu $ is
finite.

\item[(b)] For any $\lambda \in \Lambda $ and any $g\notin
H_\lambda $, the intersection $H^g_\lambda \cap H_\lambda $ is
finite.
\end{enumerate}
\end{lem}

\begin{proof}[Proof of Theorem \ref{glue}]
Obviously it suffices to deal with the case $m=1$. The general
case will follow if we apply the theorem for $m=1$ several times.

Let $t\in G$ be an arbitrary element. Passing to a new relative
generating set $X^\prime =X\cup \{ t\} $ if necessary, we may
assume that $t\in X$. By Lemma \ref{non-com} there are
non--commensurable elements $h_1,h_2\in H^0$ such that $E_G(h_1)$
and $E_G(h_2)$ are cyclic. According to Theorem \ref{E(g)}, $G$ is
hyperbolic relative to $\Hl \cup E_G(h_1)\cup E_G(h_2)$.

Let $\mu ,\e , \rho $ be
constants such that the conclusions of Lemma \ref{G1} and Lemma
\ref{torsion} are satisfied for $\lambda =1/3 $, $c=2$, $N=1$. By Theorem
\ref{w}, there are $n$ and $m_1, \ldots , m_n$ such that
the set $\mathcal R$ of all cyclic shifts and their inverses of the  word
$$R\equiv th_1^{m_1}h_2^{m_1}\ldots h_1^{m_n}h_2^{m_n}$$
in the alphabet $\mathcal A= X\cup\mathcal H\cup
(E_G(h_1)\setminus \{ 1\} )\cup (E_G(h_2)\setminus \{ 1\} )$ satisfies the $C(1/3, 2, \e , \mu , \rho )$--condition (here $h_i^{m_j}$ is regarded as a letter in $E_G(h_i)\setminus \{
1\}$, $i=1,2$, $j=1, \ldots , n$).    Indeed it suffices to chose large enough $n$ and $m_1, \ldots , m_n$ satisfying $m_i\ne \pm m_j$ whenever $i\ne j$. Let $G_1$ be the quotient of $G$ obtained by
imposing the relation $R=1$ and $\eta $ the corresponding natural
homomorphism.

By Lemma \ref{G1}, $G_1$ is hyperbolic relative to the images of
$H_\lambda , \lambda \in \Lambda $ and $E_G(h_1)$, $E_G(h_2)$. As
any elementary group is hyperbolic, $G_1$ is also hyperbolic
relative to $\{ \eta (H_\lambda \} _{\lambda \in \Lambda }\} $
according to Lemma \ref{exhyp}. The inclusion $\eta (t)\in \eta
(H)$ follows immediately from the equality $R=1$ in $G_1$. The
third assertion of the theorem follows from Lemma \ref{G1} b) as
any element from the union $\bigcup\limits_{\lambda \in \Lambda }
H_\lambda $ has length $1$.

Similarly as $\eta $ is injective on $E_G(h_1)\cup E_G(h_2)$,
$\eta (h_1)$ and $\eta (h_2) $ are elements of infinite order.
Note also that $\eta (h_1)$ and $\eta (h_2) $ are not
commensurable in $G_1$. Indeed otherwise the intersection $\big(
\eta (E_G(h_1))\big) ^g\cap \eta (E_G(h_2))$ is infinite for some
$g\in G$ contradictory the first assertion of Lemma \ref{malnorm}.
Assume now that $g\in E_{G_1}(\eta (h_1))$, where $E_{G_1}(\eta
(h_1))$ is the maximal elementary subgroup of $G_1$ containing
$\eta (h_1)$. By the first assertion of Theorem \ref{E(g)}, $\big(
\eta (h_1^m)\big) ^g= \eta (h_1^{\pm m})$ for a certain $m\ne 0$.
Therefore, $\big( \eta (E_G(h_1)) \big) ^g\cap \eta (E_G(h_1))$
contains $\eta (h_1^m)$ and in particular this intersection is
infinite. By the second assertion of Lemma \ref{malnorm}, this
means that $g\in \eta (E_G(h_1))$. Thus we proved that $\eta
(E_G(h_1))=E_{G_1}(\eta (h_1))$. The same is true for $h_2$.
Finally, using injectivity of $\eta $ on $E_G(h_1)\cup E_G(h_2)$
again, we obtain $$E_{G_1}(\eta (h_1))\cap E_{G_1}(\eta
(h_2))=\eta (E_G(h_1))\cap \eta (E_G(h_2))=\eta \big( E_G(h_1)\cap
E_G(h_2)\big) = \{ 1\} .$$ This means that the image of $H$ is a
suitable subgroup of $G_1$. To complete the proof it remains to
note that the last assertion of the theorem follows from Lemma
\ref{torsion}.
\end{proof}

%%%%%%%%%%%%%%%%%%%%%%%%%%%%%%%%%%%%%%%%%%%%%%%%%%%%%%%%%%%%%%%%

\section{Appendix. The proof of Lemma \ref{Gr0}}

%%%%%%%%%%%%%%%%%%%%%%%%%%%%%%%%%%%%%%%%%%%%%%%%%%%%%%%%%%%%%%%%

Following the referee's recommendation, we provide here the proof of Lemma \ref{Gr0}, which is a particular case of Lemma \ref{65} below. As we mentioned in Section 4, Lemma \ref{Gr0} (as well as Lemma \ref{65}) is, in fact, proved in \cite{Ols2}, although it is stated in a slightly different way there. The proof below follows \cite{Ols2} with little improvements. We stress that all results of this section should be credited  to Olshanskii.

Throughout the Appendix, let $G$ denote a group with a presentation (\ref{ZP}). Suppose that the Cayley graph $\Gamma (G, \mathcal A)$ of $G$ is hyperbolic.

We start with an auxiliary result. Recall that a graph is {\it simple} if it has no loops and multiple edges.

\begin{lem}\label{Eul}
Let $\Xi $ be a simple planar graph. Suppose that some non-negative weights $\nu (o)$, $\nu (e)$ are assigned to each vertex $o$ and each edge $e$ of $\Xi  $. Assume also that there exists a constant $a$ such that $\nu (e)\le a\nu (o)$ for any incident edge $e$ and vertex $o$. Then the sums $\nu _0$, $\nu _1$ of the weights of all vertices and edges of $\Phi $, respectively, satisfy $\nu _1\le 5a\nu _0.$
\end{lem}

\begin{proof}
By the well-known consequence of the Euler Formula, every simple planar graph contains a vertex of degree at most $5$.  The statement of the lemma easily follows from this by induction.
\end{proof}

The proof of Lemma \ref{Gr0} is divided into a sequence of lemmas. In what follows, let $G$ be a group given by (\ref{ZP}) such that the Cayley graph $\Gamma (G, \mathcal A)$ of $G$ is hyperbolic. We fix $\lambda\in (0,1]$, $c\ge 0$, $\mu \in (0, 1/16)$ and take
\begin{equation}\label{eck}
\e > c_1+2\kappa,
\end{equation}
where $c_1=c_1(\delta )$ and $\kappa=\kappa (\delta, \lambda , c)$ are the constants provided by Lemma \ref{N123} and Lemma \ref{qg} respectively. Let (\ref{quot}) be a presentation that satisfies the $\C $-condition, where $\rho $ is sufficiently large. (The exact lower bound for $\rho $ which ensures that all arguments below are correct can be easily extracted from the proofs.)

Let also $\Delta $ be a diagram over (\ref{quot}). Below we prove Lemma \ref{Gr0} by induction on the number of $\mathcal R$-cells in $\Delta $. To complete the inductive step, we will have to deal with the case when $\partial \Delta $ consists of at most $4$ quasi-geodesic segments. More precisely, let $\partial \Delta =q_1\cdots q_r$ for some $1\le r\le 4$, where $q_1, \ldots , q_r$ are $(\lambda , c)$-quasi-geodesic. We call the subpaths $q_1, \ldots, q_r$ {\it sections} of $\partial \Delta $.

\begin{defn}\label{M}
A set $\mathcal M$ of $\e $-contiguity subdiagrams of cells to cells or
cells to sections of $\partial\Delta $ is
called {\it distinguished}, if
\begin{enumerate}
\item[(a)] Distinct subdiagrams in $\mathcal M$ are disjoint.

\item[(b)] The sum of lengths of contiguity arcs of all
subdiagrams from $\mathcal M$ is not less than the same sum for
any other collections satisfying (a).

\item[(c)] $\mathcal M$ contains minimal number of subdiagrams
among all collections satisfying (a), (b).
\end{enumerate}

Further let $\Gamma \in \mathcal M$ be a $\e $-contiguity subdiagram
whose boundary is decomposed according to (\ref{pathp}) and
(\ref{sides}). The paths $s_1, s_2$ are called {\it side arcs} of $\Gamma $. If $\Gamma $ is an $\e $-contiguity subdiagram of a cell
$\Pi _1$ to a section of $\partial \Delta$,  the contiguity arc
$q_1$ of $\partial \Pi $ (and all its edges) is called {\it
outer}. If $\Gamma $ is an $\e$-contiguity subdiagram of a cell $\Pi _1$
to a cell $\Pi_2$, $q_1$ (and all its edges) is called {\it
inner}. The edges of $\Pi_1$ that are neither inner nor outer are
called {\it unbound}. Every maximal subpath of $\partial \Pi _1$
consisting of unbound edges is called an {\it unbound arc} of $\Pi
_1$. The notion of an unbound arc of a section of $\partial \Delta $ is
defined in the same way.
\end{defn}

\begin{figure}[width=100mm, hight=70mm]
\hspace{2mm}\includegraphics{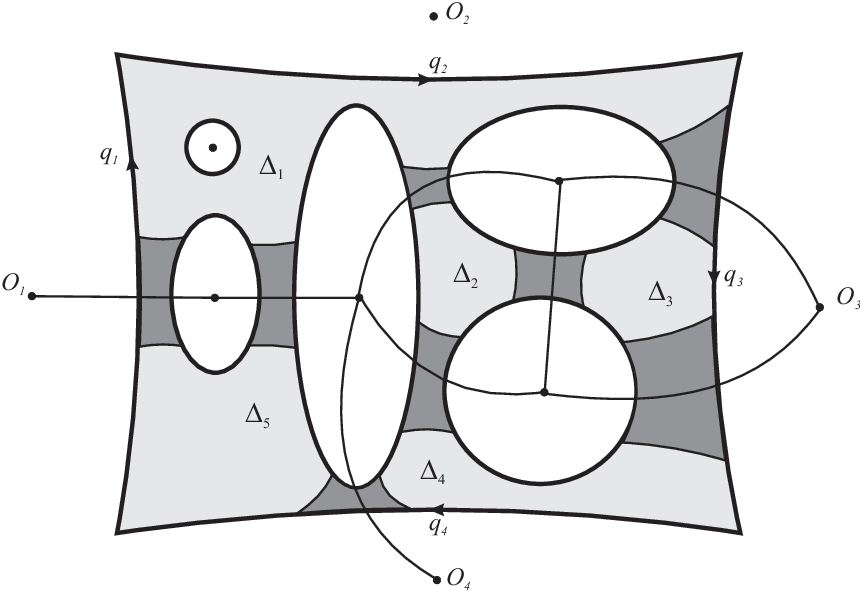}\\
\vspace{-2mm}
 \caption{A distinguished set of $\e $-contiguity subdiagrams (dark grey) in a diagram $\Delta $ with five $\mathcal R$-cells (white) and the corresponding graph $\Phi ^\prime _{\mathcal M}$.}
\label{diag}
\end{figure}

To each distinguished set $\mathcal M$ of $\e $-contiguity diagrams, we associate a planar graph $\Phi _{\mathcal M}$ as follows.
Choose a point $o(\Pi )$ inside every $\mathcal R$--cell $\Pi $ of
$\Delta $. The set of vertices of $\Phi _{\mathcal M}$ consists of all such
$o(\Pi )$. If there exists an $\e$-contiguity subdiagram $\Gamma \in
\mathcal M$ of a cell $\Pi _1$ to a cell $\Pi _2$, the two
vertices $o(\Pi _1)$ and $o(\Pi _2)$ are connected by an
(unoriented) edge through $\Gamma $. Further let $\Phi ^\prime _{\mathcal M}$ be the graph obtained from $\Phi _{\mathcal M}$ in the following way. For each $i=1, \ldots , r$, we add to $\Phi _{\mathcal M}$ a vertex $O_i$ outside of $\Delta $ and for each $\e $-contiguity subdiagram $\Gamma $ of a cell $\Pi $ to $q_i$, we connect $o(\Pi )$ to $O_i$ by an edge passing through $\Gamma $. Clearly $\Phi ^\prime _{\mathcal M}$ is also planar (see Fig. \ref{diag}).

Let now  $\Delta $ be a reduced diagram over (\ref{quot}), which has $n\ge 1$ $\mathcal R$-cells. The next $4$ results are proved under the following additional assumption. It will be eliminated later in Lemma \ref{64}.
$$
(\ast ) \;\;\; \begin{array}{l}\mbox{\it For any distinguished system of $\e$-contiguity subdiagrams $\mathcal M$ in $\Delta $, the graph $\Phi _{\mathcal M}$ is} \\ \mbox{\it  simple and inside every $2$-gon of $\Phi ^\prime _{\mathcal M}$, there is a vertex of $\Phi _{\mathcal M}$.}\end{array}
$$

Let also $\mathcal M$ be a distinguished system of $\e$-contiguity subdiagrams of $\Delta $.   Cutting off all $\mathcal R$-cells and subdiagrams $\Gamma \in \mathcal M$, we obtain a set of diagrams $\Delta _1, \ldots , \Delta _d$ over (\ref{ZP}) (see Fig. \ref{diag}). Each of them may have holes. The boundary $\partial \Delta _i$ of each diagram may be thought of as a union of $n_i$ arcs, where each arc is of one of the following types.
\begin{enumerate}
\item[(A1)] An unbound arc of an $\mathcal R$-face.

\item[(A2)] An unbound arc of $\partial \Delta $.

\item[(A3)] A side arc of some $\e $-contiguity subdiagram from $\mathcal M$.
\end{enumerate}

\begin{lem}\label{61}
Suppose that $\Delta $ satisfies ($\ast $). In the notation introduced above, we have $\sum\limits_{i=1}^d n_i\le 53 n.$
\end{lem}

\begin{proof}
Let $v,e,f$ denote the number of vertices, edges, and regions of (the planar realization of) $\Phi _{\mathcal M}^\prime $, respectively. By ($\ast $), every region of it (except possibly for the outer one) has degree at least three. Hence $f\le 2e/3+1$. By the Euler formula, we have $e\le v+f-2\le v + 2e/3-1$. This implies $e\le 3(v-1)\le 3(n+3)\le 12 n$.

Since unbound arcs of boundaries of $\mathcal R$-cells and sections of $\partial \Delta $ are separated by contiguity arcs, the total number of arcs of type (A1) or (A2) in $\Delta $ is not greater than $$2|\mathcal M|+n+r\le 2e +n+4\le 24 n +n+4 \le 29n.$$ Clearly the number of arcs of type (A3) is at most $2|\mathcal M|=2e \le 24n$. So the total number of arcs of types (A1)-(A3) is at most $53n$.
\end{proof}

All results below are proved up to passing from $\Delta $ to an $\mathcal O$-equivalent diagram (see the definition before Lemma \ref{Gr0}).

\begin{lem}\label{62}
Suppose that $\Delta $ satisfies ($\ast $). Let $S$ denote the sum of lengths of all unbound arcs of type (A1) in $\Delta $. Then $S<n\sqrt{\rho }$.
\end{lem}

\begin{proof}
Let $S_i$  denote the sum of lengths of all arcs of type (A1) in $\partial \Delta _i$, $i=1, \ldots , d$. Assume that $S\ge n\sqrt{\rho }$. Then
\begin{equation}\label{621}
S_i\ge n_i\sqrt{\rho}/60
\end{equation}
 for some $i$. Indeed otherwise we have
$$
S=\sum\limits_{i=1}^d S_i\le\frac{\sqrt{\rho}}{60} \sum\limits_{i=1}^d n_i
< n\sqrt{\rho}
$$
By Lemma \ref{61}.

For every $i=1, \ldots , d$, the diagram $\Delta _i$ has at most $n_i$ holes. Hence we can cut it by $l\le n_i$ paths $t_1, \ldots, t_l$ into a simply connected diagram $\widetilde{\Delta _i}$. We assume that the collection of cutting paths $t_1, \ldots , t_l$ is chosen so that the sum $\sum\limits_{i=1}^l l(t_i) $ is minimal. The boundary of $\widetilde{\Delta _i}$ decomposes into $k_i$ subpaths, each of which is either $t_j^{\pm 1}$ for some $j=1, \ldots , l$, or an arc of one of the types (A1)-(A3) (or a part of such an arc arising after cutting along $t_1, \ldots , t_l$). Therefore,
\begin{equation}\label{622}
k_i\le 4n_i.
\end{equation}

Note also that passing to a diagram which is $\mathcal O$-equivalent to $\Delta $ if necessary, we may assume that the labels of $t_i$'s are geodesic in $G$ without loss of generality. Indeed if $w_i$ is a geodesic word representing the same element as $\phi (t_i)$ in $G$, let $\Sigma _i$ be a diagram over (\ref{ZP}) with boundary label $w_i(\phi (t_i))^{-1}$. Further let $\Xi _i$ be the diagram obtained by gluing $\Sigma _i$ and its mirror copy along $w_i$. Clearly $\phi \partial (\Xi_i)=\phi (t_i)(\phi (t_i))^{-1}$. We can use $0$-refinement (see the subsection ``Van Kampen diagrams" in Section 3) to create a copy $t_i^\prime $ of the paths $t_i$ in $\Delta $ with the same label and endpoints such that $t_i(t_i^\prime )^{-1}$ bounds a subdiagram consisting of $0$-cells. Further we cut this subdiagram and fill in the obtained hole with $\Xi_i$. After this transformation, the vertices $(t_i)_-$ and $(t_i)_+$ are connected by a path whose label $w_i$ is geodesic in $G$, and we can replace $t_i$ with that paths. Note that this transformation does not affect $\partial \Delta $, $\mathcal R$-cells, and distinguished $\e$-contiguity subdiagrams.

Similarly we may assume that for every (sub)arc $p$ of an arc of type (A1)-(A3) in $\partial \widetilde{\Delta _i}$, the diagram $\widetilde{\Delta _i}$ contains a path $\bar p$ with the same endpoints as $p$ such that $\phi (\bar p)$ is geodesic in $G$ and $p(\bar p)^{-1}$ bounds a subdiagram over (\ref{ZP}). Since $\widetilde{\Delta _i}$ is a simply connected diagram over (\ref{ZP}), its $1$-skeleton can be naturally mapped to $\Ga$ by a map preserving labels and orientation. In what follows we keep the same notation $t_i^{\pm 1}$, $p$, etc., for the images of the paths $t_i^{\pm 1}$, $p$, etc., in $\Ga $ (although $t_i^{-1}$ is not the inverse path of $t_i$ there). In particular, $p$ and $\bar p$ belong to the closed $\kappa $-neighborhoods of each other in $\Ga$ by Lemma \ref{qg}.

In the notation of Lemma \ref{N123}, let $N_2$ consist of images of segments $t_i^{\pm 1}$, $i=1, \ldots , l$, in $\Ga $ and images of $\bar p$ for (sub)arcs $p$ of type (A2). Further let $N_1$ and $N_3$ consist of images of $\bar p$ for (sub)arcs $p$ of type (A1) and (A3), respectively.

Clearly $\sigma _3\le \e k_i$. On the other hand, by (\ref{621}) and (\ref{622}) we obtain
$$
\begin{array}{rl}
\sigma _1 & = \sum \limits_{\bar p\in N_1} l(p_1) \ge \sum \limits_{\bar p\in N_1} (\lambda l(p) - c)\ge \lambda S_i-k_ic> \lambda n_i \sqrt{\rho}/60 -k_ic \\ & \\ & \ge
\lambda k_i\sqrt{\rho} /240 -k_i c =k_i(\lambda \sqrt{\rho } /240 - c).
\end{array}
$$
Taking large enough $\rho $, we may assume that $a=\lambda \sqrt{\rho } /240 - c>\max\{ 1000\e , c_2\} $. Hence by Lemma \ref{N123}, there exist subsegments $\bar q_j$, $j=1,2$, of  some $\bar p_1\in N_1$ and $\bar p_2\in N_1\cup N_2$, respectively, such that $l(p_j)\ge 10^{-3} a$, $j=1,2$, and the distance between the respective endpoints of $p_1$ and $p_2$ is at most $c_1$. This means that there exist
subsegments $q_j$ of $p_j$, $j=1,2$, of lengths at least  $10^{-3}a -2\kappa$ such that
\begin{equation}\label{q1q2dist}
\max\{ \dxh ((q_1)_-, (q_2)_-),\, \dxh ((q_1)_+, (q_2)_+)\} \le 2\kappa +c_1<\e.
\end{equation}
The later inequality uses (\ref{eck}). Now there are two cases to consider.

{\it Case 1.} $\bar p_2\in N_1$. That is, $q_1$ and $q_2$ are subpaths of boundaries of some $\mathcal R$-cells in $\Delta $. After passing to an $\mathcal O$-equivalent diagram if necessary, we may assume that there exists a subdiagram $\Gamma $ of $\Delta $ such that $\partial \Gamma =q_1s_1q_2^{-1}s_2$, where $l(s_j)<\e$, $j=1,2$, and $\Gamma $ does not intersect $\mathcal R$-cells and distinguished $\e $-contiguity subdiagrams. This contradicts the maximality of $\mathcal M$.

\begin{figure}
   \hspace{27mm} \includegraphics{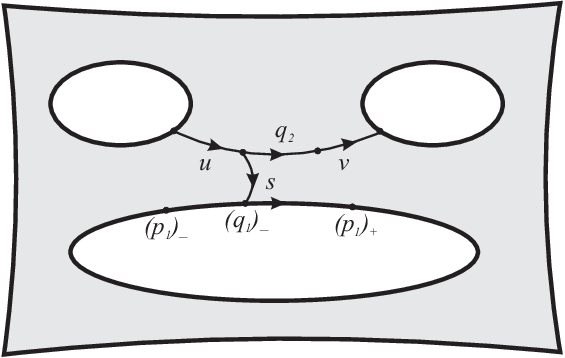}\\
   \vspace{-3mm}
  \caption{Replacing a cutting path with a shorter one.}\label{cut}
\end{figure}

{\it Case 2.} $\bar p_2\in N_2$. If $q_2$ is a subpaths of $\partial \Delta $, we get a contradiction as above by constructing a new $\e $-contiguity subdiagram of an $\mathcal R$--cell to $\partial \Delta $. Thus we only need to deal with the case $p_2=t_j^{\pm 1}$. Without loss of generality, we may assume that $p_2=t_j$. Let $t_j=uq_2v$. By (\ref{q1q2dist}), after passing to an $\mathcal O$-equivalent diagram if necessary, we can find a paths $s$ of lengths at most $\e$ in $\Delta $ connecting $(q_2)_-$ to $(q_1)_-$ (see Fig. \ref{cut}). Note that
$$l(us)\le l(t_j)- l(q_2) +\e \le l(t_j) - (10^{-3}(\lambda \sqrt{\rho } /240 - c) -2\kappa)  +\e <l(t_j)$$ if $\rho $ is large enough. This contradicts our assumption that
$\sum\limits_{i=1}^l l(t_i) $ is minimal.

Thus in both cases we obtain a contradiction. Hence $S< n\sqrt{\rho }$.
\end{proof}

\begin{lem}\label{63}
Suppose that $\Delta $ satisfies ($\ast $). Let $\Sigma _0$ and $\Sigma $ denote the sum  of lengths of all outer arcs and the sum of perimeters of all $\mathcal R$-cells in $\Delta $, respectively. Then
$\Sigma _0 > (1-11\mu )\Sigma $.
\end{lem}

\begin{proof}
For each $\mathcal R$-cell $\Pi $ of $\Delta $, we assign the weight $\nu (o)=l(\partial \Pi )$ to the corresponding vertex $o$ of the estimating graph $\Phi _\mathcal M$. Clearly the sum of weights of all vertices in $\Phi _\mathcal M$ equals $\Sigma $. Further let $\Gamma \in \mathcal M$  be a  distinguished $\e $-contiguity subdiagram with $\partial \Gamma =s_1q_1s_2q_2$, where $q_1, q_2$ are the contiguity arcs. We assign the weight $\nu (e)= l(q_1)+l(q_2)$ to the edge $e$ corresponding to $\Gamma $. Note that the sum of weights of all edges in $\Gamma \in \mathcal M$ is equal to the sum of lengths $\Sigma _{inn}$ of all inner arcs in $\Delta $.

By Lemma \ref{O52}, $\nu (e)\le 2\mu \nu (o)$ whenever $o$ and $e$ are incident. Hence $\Sigma _{inn}\le 10\mu\Sigma $ by Lemma \ref{Eul} and ($\ast $). By Lemma \ref{62}, we have $S<\Sigma /\sqrt{\rho }$, where $S$ is the sum of lengths of all unbound arcs, since the perimeter of every $\mathcal R$-cell in $\Delta $ is at least $\rho $ by the $\C $-condition. Therefore, $$\Sigma _0=\Sigma -\Sigma _{int} - S>\Sigma - 10\mu\Sigma -\Sigma/\sqrt{\rho } >(1-11\mu )\Sigma $$ if $\rho >\mu ^{-2}.$
\end{proof}

The following corollary is immediate.

\begin{cor}\label{64}
Suppose that $\Delta $ satisfies ($\ast $). Passing to an $\mathcal O$-equivalent diagram if necessary, we can find an $\mathcal R$-cell $\Pi $ of $\Delta $ and disjoint $\e $-contiguity subdiagrams $\Gamma _{i,j}$ of $\Pi $ to sections $q_j$, $j=1, \ldots , r$ of $\partial \Delta $ such that $\sum\limits_{i,j} (\Pi, \Gamma _{i,j}, q_j)>1-11\mu $.
\end{cor}

We are now ready to eliminate assumption ($\ast $) and prove the main result of the Appendix, of which Lemma \ref{Gr0} is a particular case corresponding to $r=1$. For convenience, we recall all assumptions here.

\begin{lem}\label{65}
Let $G$ be a group with a presentation (\ref{ZP}). Suppose that the Cayley graph $\Gamma (G, \mathcal A)$ of $G$ is hyperbolic. Then for any $\lambda\in (0,1]$, $c\ge 0$, and $\mu \in (0, 1/16]$, there exist
$\e \ge 0$ and $\rho >0$ with the following property.
Let $\mathcal R$ be a symmetrized set of words in $\mathcal A$
satisfying the $C(\e, \mu , \lambda , c , \rho )$--condition,
$\Delta $ a reduced van Kampen diagram over the presentation
(\ref{quot}) such that $\partial \Delta =q_1\cdots q_r$ for some $1\le r\le 4$, where $q_1, \ldots , q_r$ are $(\lambda , c)$-quasi-geodesic. Assume that $\Delta $ has at least one $\mathcal R$--cell. Then up to passing to an $\mathcal O$-equivalent diagram, the following conditions hold.
\begin{enumerate}
\item[(a)] The diagram $\Delta $ satisfies ($\ast $).

\item[(b)] There is an $\mathcal R$-cell $\Pi $ of $\Delta $ and disjoint $\e $-contiguity subdiagrams $\Gamma _{j}$ of $\Pi $ to sections $q_j$, $j=1, \ldots , r$, of $\partial \Delta $ (some of them may be absent) such that $\sum\limits_{j=1}^r (\Pi, \Gamma _{j}, q_j)>1-13\mu $.
\end{enumerate}
\end{lem}

\begin{proof}
We proceed by induction on $n$, the number of $\mathcal R$-cells in $\Delta $.

(a) Suppose that there are multiple edges in $\Phi _{\mathcal M}$. That is, there are two distinguished $\e $-contiguity subdiagrams $\Theta _1$ and $\Theta _2$ between some $\mathcal R$-cells $\Pi _1$ and $\Pi _2$. Consider a subdiagram $\Xi $ in $\Delta $ such that: (i) $\partial \Xi =s_1t_1s_2t_2$, where $s_j$ is a side arc of $\Theta _j$ and $t_j$ is a subpaths of $\partial \Pi _j$, $j=1,2$; (ii) $\Xi $ does not contain $\Pi _1$ and $\Pi _2$. By the definition of $\mathcal M$, we can not include $\Theta _1$ and $\Theta _2$ into a single $\e $-contiguity subdiagram of $\Pi _1$ to $\Pi _2$. This means that $\Xi$ contains at least one $\mathcal R$-cell. On the other hand, the number of $\mathcal R$-cells in $\Xi $ is smaller than $n$ according to (ii).  Hence by the inductive assumption,  there is an $\mathcal R$-cell $\Pi $ in $\Xi $ and $\e $-contiguity subdiagrams $\Gamma _1, \ldots , \Gamma _4$ of $\Pi $ to $s_1, t_1, s_2, t_2$, respectively, such that $$(\Pi ,\Gamma _1, s_1)+\cdots +(\Pi ,\Gamma _4, t_2)> 1-13\mu . $$

On the other hand $(\Pi ,\Gamma _2, t_1)+(\Pi ,\Gamma _4, t_2)<2\mu $ by Lemma \ref{O52} since $\Gamma _2$ and $\Gamma _4$ are $\e $-contiguity subdiagrams of $\Pi $ to $\Pi _1$ and $\Pi _2$, respectively. Further if $u$ is a contiguity arc of $\Pi $ to $s_1$, then $l(u)\le \lambda ^{-1} (3\e +c)$ since $u$ is $(\lambda , c)$ quasi-geodesic and $l(s_1)\le \e$. Therefore, $(\Pi, \Gamma _1, s_1)=l(u)/l(\partial \Pi )<l(u)/\rho < \mu/2 $ if $\rho $ is large enough. Similarly  $(\Pi, \Gamma _3, s_2)<\mu /2$. Thus $$(\Pi ,\Gamma _1, s_1)+\cdots +(\Pi ,\Gamma _4, t_2)<3\mu . $$ We obtain a contradiction since $1-13\mu >3\mu$ for $\mu < 1/16 $.

Thus the graph $\Phi _{\mathcal M}$ does not have multiple edges. Arguing as above, it is easy to show that $\Phi _{\mathcal M}$ can not contain loops either. The only difference is that the subdiagram $\Xi $ will be bounded by $st$, where $l(s)<\e $ and $t$ is a subpaths of $\partial \Pi_1$. (Here $\Pi _1$ is the $\mathcal R$-cell of $\Delta $ such that there is a lop incident to the corresponding vertex in $\Phi _{\mathcal M}$.)

Finally inside every $2$-gon $ef$ of $\Phi _{\mathcal M}^\prime $ there is a vertex of $\Phi _{\mathcal M}$ since otherwise one can include the $\e$-contiguity subdiagrams corresponding to the edges $e$ and $f$ of $\Phi _{\mathcal M}^\prime $ into a single $\e$-contiguity subdiagram in the obvious way, contrary to the definition of $\mathcal M$.

(b) By (a), we can choose an $\mathcal R$-cell $\Pi $ in $\Delta $ (passing to an $\mathcal O$-equivalent diagram if necessary) and subdiagrams $\Gamma _{i,j}$ satisfying the inequality in Corollary \ref{64}. Let us consider the subdiagram $\Gamma ^1$ of $\Delta $ such that: (i) $\partial \Gamma ^1=s_1t_1s_2t_2$, where $s_1, s_2$ are side arcs of some of $\Gamma_{i,1}$'s and $t_1, t_2$ are subpaths of $\partial \Pi $ and the section $q_1$ of $\partial \Delta $, respectively; (ii) $\Gamma ^1$ contains all $\Gamma _{i, 1}$'s. Let $m_1$ be the number of $\mathcal R$-cells in $\Gamma ^1$. We similarly construct $\Gamma^2, \ldots , \Gamma ^r$ and define $m_2, \ldots , m_r$.

Suppose that the cell $\Pi $ is chosen in such a way that the sum $m(\Pi )=m_1+\ldots + m_r$ is minimal among all cells satisfying the condition  $\sum\limits_{i,j} (\Pi, \Gamma _{i,j}, q_j)>1-11\mu $. If $m_1=\cdots =m_r=0$, then each of the sets $\{ \Gamma _i,1\}$, $\ldots $, $\{ \Gamma _i,1\}$ consists of at most one diagram. Indeed otherwise one could include at least $2$ of them into a single $\e$-contiguity subdiagram, which contradicts the definition of $\mathcal M$.

Thus we may assume that $m_1>0$. Then by Corollary \ref{64} and the inductive assumption, the subdiagram $\Gamma ^1 $ contains an $\mathcal R$-cell $\Pi ^\prime  $ and $\e $-contiguity subdiagrams $\Gamma ^\prime_{i,1} , \ldots , \Gamma ^\prime  _{i, 4}$ of $\Pi^\prime $ to $s_1$, $t_1$, $s_2$, $t_2$, respectively, such that
\begin{equation}\label{651}
\sum\limits_{i=1}^{k_1} (\Pi^\prime, \Gamma ^\prime_{i,1}, s_1)+ \cdots +\sum\limits_{i=1}^{k_4} (\Pi^\prime, \Gamma ^\prime_{i,4}, t_2) >1-11\mu
\end{equation}
Note that $k_2\le 1$ since all $\Gamma ^\prime_{i,2}$'s are $\e$-contiguity subdiagrams of $\Pi ^\prime  $ to $\Pi $ and $\Delta $ satisfies ($\ast $) by (a). If $k_2=1$, we have $(\Pi^\prime, \Gamma ^\prime  _{1,2}, t_1)<\mu $ by Lemma \ref{O52}. Further if $k_1\ge 1$, we can construct an $\e$-contiguity subdiagram $\Gamma $ of $\Pi ^\prime  $ to $t_1$ such that: (i) $\partial \Gamma =st$, where $l(s)\le 3\e $ and $t$ is a subpaths of $\partial \Pi ^\prime  $; (ii) $\Gamma $ contains at least one $\mathcal R$-cell. This leads to a contradiction as in part (a). Thus $k_1\le 1$ and if $k_1=1$ we obtain $(\Pi ^\prime, \Gamma ^\prime  _{i,1}, s_1)<\mu /2$ as in part (a). Similarly $k_3\le 1$ and if $k_3=1$, we have $(\Pi^\prime, \Gamma ^\prime  _{i,3}, s_2)<\mu /2$. These inequalities together with (\ref{651}) imply $\sum\limits_{i=1}^{k_4} (\Pi^\prime, \Gamma ^\prime _{i,4}, t_2) >1-13\mu $. However $m(\Pi ^\prime  )<m(\Pi )$ since the cell $\Pi ^\prime  $ counts in $m (\Pi )$ but not in $m(\Pi ^\prime  )$. This contradicts the choice  of $\Pi $.
\end{proof}

\vspace{1cm}

\noindent \small  Stevenson Center 1326, Department of
Mathematics, Vanderbilt University\\  Nashville, TN 37240, USA

\vspace{3mm}

\noindent \small \it E-mail address: \tt
denis.osin@gmail.com
\end{document}